\documentclass{amsart}
\usepackage{amscd,amsfonts,amssymb}
\usepackage{textgreek}
\usepackage{tikz}
\usepackage{pgfplots}
\usetikzlibrary{automata,positioning,calc,trees}
\usetikzlibrary{intersections,pgfplots.fillbetween}
\usepackage{graphicx,tikz}
\usepackage{pgf,tikz,pgfplots}
\usepackage{mathrsfs}
\usetikzlibrary{arrows}
\usepackage{enumerate}
\usepackage[shortlabels]{enumitem}
\usepackage{mathrsfs}
\usepackage{amssymb,amsmath,amsthm,color}
\usepackage{caption,subcaption}
\usepackage{hyperref}
\usepackage{cleveref}
\usepackage[alphabetic,bibtex-style]{amsrefs}
\usepackage{url}
\usepackage{setspace}
\usepackage{float}
\usepackage{makecell}
\tikzstyle{startstop} = [rectangle, rounded corners, minimum width=3cm, minimum height=1cm,text centered, draw=black]
\tikzstyle{process} = [rectangle, rounded corners, minimum width=3cm, minimum height=1cm, text centered, draw=black]
\tikzstyle{arrow} = [thick,->,>=stealth]

\BibSpec{article}{%
	+{}  {\PrintAuthors}                {author}
	+{,} { \textit}                     {title}
	+{.} { }                            {part}
	+{:} { \textit}                     {subtitle}
	+{,} { \PrintContributions}         {contribution}
	+{.} { \PrintPartials}              {partial}
	+{,} { }                            {journal}
	+{}  { \textbf}                     {volume}
	+{}  { \PrintDatePV}                {date}
	+{,} { \issuetext}                  {number}
	+{,} { \eprintpages}                {pages}
	+{,} { }                            {status}
	+{,} { \url}                        {url}    % <---- ADDED
	+{,} { \PrintDOI}                   {doi}
	+{,} { available at \eprint}        {eprint}
	+{}  { \parenthesize}               {language}
	+{}  { \PrintTranslation}           {translation}
	+{;} { \PrintReprint}               {reprint}
	+{.} { }                            {note}
	+{.} {}                             {transition}
}

\date{\today}
\allowdisplaybreaks
\newcommand{\R}{{\mathbb R}}       % Field of real numbers

\newcommand{\Z}{{\mathbb Z}}       % Ring of integer numbers

\newtheorem{theorem}{Theorem}[section]
\newtheorem{lemma}[theorem]{Lemma}

\newtheorem{proposition}[theorem]{Proposition}

\newtheorem*{theorem*}{Theorem}

\theoremstyle{definition}

\numberwithin{equation}{section}

\allowdisplaybreaks

\begin{document}
	
	\title[Bilinear Bochner-Riesz operator at the critical index]{Weak--type $(1,1;\frac{1}{2})$--estimate for the bilinear Bochner--Riesz mean at the critical index}
	
	\author[S.S. Choudhary]{Surjeet Singh Choudhary}
	\address{Surjeet Singh Choudhary\\
		National Center for Theoretical Sciences\\
		National Taiwan University\\
		Taipei-106, Taiwan.}
	\email{surjeet19@ncts.ntu.edu.tw}
	
	\author[S. Shrivastava]{Saurabh Shrivastava}
	\address{Saurabh Shrivastava\\
		Department of Mathematics\\
		Indian Institute of Science Education and Research Bhopal\\
		Bhopal-462066, India.}
	\email{saurabhk@iiserb.ac.in}
	
	\begin{abstract}
		We prove the weak-type estimate at the end-point $\big(1,1,\frac{1}{2}\big)$ for the bilinear Bochner--Riesz mean at the critical index $\alpha_c=n-\frac{1}{2}$ in all dimensions $n\geq 1$. Our proof is motivated by Christ's method for the linear Bochner--Riesz operator at the critical index. We develop a bilinear analogue of the $TT^*$ method to address the end-point boundedness of the bilinear Bochner--Riesz operator. This approach provides a potentially useful tool in addressing the end-point issues for bilinear operators.
	\end{abstract}

	\maketitle
	\tableofcontents
	\section{Introduction}
    \subsection{Linear Bochner--Riesz means}
	Let $n\geq2$ and $\alpha\geq0$. For a Schwartz function $f$ on $\R^n$, the Bochner--Riesz operator of order $\alpha$ is defined by
	\[B^\alpha f(x)=\int_{\R^n}(1-|\xi|^2)_+^\alpha\widehat f(\xi)e^{2\pi i x\cdot\xi}\,d\xi,\]
	where $r_+=\max\{r,0\}$ and $\widehat f(\xi)=\int_{\R^n}f(x)e^{-2\pi i x\cdot\xi}\,dx$. The operators $B^\alpha$ were introduced as a summability method for multiple Fourier integrals; the parameter $\alpha$ measures how smoothly the multiplier vanishes at the unit sphere. The study of the $L^p$ boundedness of $B^\alpha$, and its connection with the restriction and Kakeya problem, is a central theme of research in harmonic analysis; see \cite{Stein93}.
	
	We can rewrite the Bochner--Riesz operator $B^\alpha$ as a convolution operator,
    \[B^\alpha f(x)=\int_{\R^n}K^\alpha(x-y)f(y)\,dy,\]
	where
	\[K^\alpha(y)=\frac{\Gamma(\alpha+1)}{\pi^\alpha}\frac{J_{\frac n2+\alpha}(2\pi|y|)}{|y|^{\frac n2+\alpha}},\]
	and $J_\nu$ is the Bessel function of order $\nu$. By the asymptotic expansion of the Bessel function, we have
	\[|K^\alpha(y)|\lesssim(1+|y|)^{-\frac{n+1}{2}-\alpha}.\]
	Thus, $K^\alpha\in L^1(\R^n)$ when $\alpha>\frac{n-1}{2}$. Consequently, $B^\alpha$ is bounded on $L^p(\R^n)$ for every $1\leq p\leq\infty$ for the range $\alpha>\frac{n-1}{2}$. For this reason, the index $\alpha_c=\frac{n-1}{2}$ is called the critical index. At the critical index, the kernel decays exactly like $|y|^{-n}$ and is no longer integrable, and $B^{\frac{n-1}{2}}$ is unbounded on $L^1$ and $L^\infty$; see \cite{Stein93, Tao99}.
	
	Below the critical index, the range of $p$ for $L^p$-boundedness of $B^\alpha$ depends on $\alpha$. For $n\geq2$ and $1\leq p\leq\infty$, $p\neq2$, the Bochner--Riesz conjecture states that $B^\alpha$ is bounded on $L^p(\R^n)$ if and only if
	\[\alpha>\alpha(p):=\max\left\{n\left|\frac{1}{p}-\frac{1}{2}\right|-\frac{1}{2},0\right\}.\]
	
	The necessity of the condition $\alpha>\alpha(p)$ is classical, so the content of the conjecture is sufficiency. When $\alpha=0$, $B^0$ is the ball multiplier operator. Fefferman \cite{Fefferman71} proved that the ball multiplier operator $B^0$ is unbounded on $L^p(\R^n)$ for every $p\neq 2$ when $n\geq2$. Carleson and Sj\"olin \cite{CS72} established the Bochner--Riesz conjecture in dimension two. In higher dimensions, the conjecture remains open, although substantial progress has been made over the years; we refer to \cite{Tomas75, Lee04, BG11, GOWWZ25} and the references therein. 
	
	At the endpoint $\alpha=\alpha(p)$, the strong-type $L^p$ estimate fails for $p\ne2$ when $n\geq2$, and it is natural to expect weak-type substitutes at the endpoint. Indeed, the weak-type estimates at this index have been studied extensively; see, for example, \cite{Christ87, Seeger962, Tao98}. At the critical index $\alpha=\frac{n-1}{2}$, which corresponds to $p=1$, Christ \cite{Christ88} proved the weak-type $(1,1)$ estimate
	\[\|B^{\frac{n-1}{2}}f\|_{L^{1,\infty}(\R^n)}\lesssim\|f\|_{L^1(\R^n)}.\]
	Vargas \cite{Vargas96} extended this result to weighted weak-type $(1,1)$ estimates for $A_1$ weights. 
    
    Christ's theorem does not follow from standard Calder\'on--Zygmund theory. At the critical index the kernel has the size $|x|^{-n}$ of a Calder\'on--Zygmund kernel. Because of its oscillation, however, its gradient decays only like $|x|^{-n}$ instead of $|x|^{-n-1}$, so H\"ormander's smoothness condition fails. In order to overcome this difficulty, he decomposed the kernel dyadically, $K^{\frac{n-1}{2}}=\sum_{j\geq0}K_j$ with $K_j$ essentially supported where $|x|\sim2^j$, and applied the Calder\'on--Zygmund decomposition $f=g+\sum_Q b_Q$. For each $j$ and each $s\geq0$, he organized the bad functions supported on cubes of side length $2^{j-s}$. Then, the main step is an $L^2$ estimate for $\sum_j K_j*\big(\sum_{\ell(Q)=2^{j-s}}b_Q\big)$ with a decay of the order $2^{-\epsilon s}$. In order to prove this estimate, he expands the square, which reduces matters to study the terms $K_j*\widetilde{K_k}$, where $\widetilde K(x)=\overline{K(-x)}$, and estimates these kernels using the oscillation of $K_j$ together with the cancellation of the functions $b_Q$. This is a $TT^*$ argument in the theory of linear operators, and it is the point of view we develop adapted to the bilinear setting. 
	
 In the same paper, Christ \cite{Christ88} proved the weak-type $(1,1)$ estimates for maximal averages over balls in dimension two with angular density $|\Omega|$, where $\Omega\in L^q(\mathbb S^1)$ for some $q>1$. Christ and Rubio de Francia \cite{CR88} extended this result to all dimensions with $\Omega\in L\log L$. When $\Omega$ also has mean zero, they established the weak-type $(1,1)$ estimates for rough singular integral operators in dimension two; this was also obtained independently by Hofmann \cite{Hofmann88}. Seeger \cite{Seeger96} later extended the result to all dimensions $n\geq2$; see also Tao \cite{Tao98}. 
 
	These weak-type estimates have since been strengthened through the sparse domination principle.  Conde-Alonso, Culiuc, Di Plinio and Ou \cite{CCDO17} proved sparse bounds both for the Bochner--Riesz operator at the critical index and for rough singular integral operators. Their proof combines a general sparse domination principle with the oscillatory estimates from Christ's argument for the Bochner--Riesz operator and with estimates from Seeger's work on rough singular integrals. Sparse bounds for Bochner--Riesz multipliers below the critical index were studied by Lacey, Mena and Reguera \cite{LMR19}, see also \cite{BeneaBernicot}. In the plane, Kesler and Lacey \cite{KL18} obtained endpoint sparse bounds for Bochner--Riesz multipliers of order $\alpha(p)$. More recently, Beltran, Roos and Seeger \cite{BRS25} obtained endpoint sparse bounds that are almost optimal in the plane, together with partial results in higher dimensions.

	\subsection{Bilinear Bochner-Riesz means}
	Let $n\geq1$. The bilinear Bochner--Riesz operator of order $\alpha\geq0$ is defined by 
	\[\mathcal B^\alpha(f,g)(x)=\int_{\R^n}\int_{\R^n}(1-|\xi|^2-|\eta|^2)_+^\alpha\widehat f(\xi)\widehat g(\eta)e^{2\pi i x\cdot(\xi+\eta)}\,d\xi\,d\eta,\]
	where $f,g\in\mathcal S(\R^n)$. 
    
    Its multiplier is the linear Bochner--Riesz multiplier on $\R^{2n}$, evaluated at $(\xi,\eta)$. In kernel form, we can express $\mathcal B^\alpha$ as follows 
	\[\mathcal B^\alpha(f,g)(x)=\int_{\R^n}\int_{\R^n}K^\alpha(y,z)f(x-y)g(x-z)\,dy\,dz,\]
	where
	\[K^\alpha(y,z)=\frac{\Gamma(\alpha+1)}{\pi^\alpha}\frac{J_{n+\alpha}(2\pi|(y,z)|)}{|(y,z)|^{n+\alpha}}.\]
	As in the linear case, the asymptotic expansion of the Bessel function gives
	\[|K^\alpha(y,z)|\lesssim(1+|(y,z)|)^{-n-\alpha-\frac{1}{2}}.\]
	When $\alpha>n-\frac{1}{2}$, the right-hand side in the above is bounded by $(1+|(y,z)|)^{-2n-2\epsilon}$ for some $\epsilon>0$, which is at most $(1+|y|)^{-n-\epsilon}(1+|z|)^{-n-\epsilon}$. Thus $|\mathcal B^\alpha(f,g)|$ is dominated pointwise by the product of two convolution operators with integrable kernels on $\R^n$. Consequently, we get 
	\[\|\mathcal B^\alpha(f,g)\|_{L^p(\R^n)} \lesssim\|f\|_{L^{p_1}(\R^n)}\|g\|_{L^{p_2}(\R^n)}\]
	for all $1\leq p_1,p_2\leq\infty$ satisfying $\frac 1p=\frac1{p_1}+\frac1{p_2}$. Accordingly, $\alpha_c=n-\frac{1}{2}$ is the critical index for the bilinear Bochner--Riesz operator. It coincides with the critical index $\frac{2n-1}{2}$ of the linear Bochner--Riesz operator on $\R^{2n}$.
	
	When $\alpha=0$, the operator $\mathcal B^0$ is the bilinear ball multiplier operator. Grafakos and Li \cite{GL06} proved that the bilinear ball multiplier is bounded in dimension one in the local $L^2$ range $2\leq p_1,p_2,p'<\infty$, where $p'$ denotes the conjugate exponent of $p$ given by $1=\frac1{p}+\frac1{p'}$. In dimensions $n\geq2$, Diestel and Grafakos \cite{DG07} showed that $\mathcal B^0$ is unbounded in the Banach range whenever exactly one of $p_1,p_2,p'$ is less than $2$. These results parallel Fefferman's theorem in the linear case. Bernicot, Grafakos, Song and Yan \cite{BGSY15} initiated the study of $\mathcal B^\alpha$ for positive $\alpha$. For a fixed $\alpha$, the range of boundedness was later improved by Jeong, Lee and Vargas \cite{JLV18} in the Banach range and by Liu and Wang \cite{LW20} in the non-Banach range $p<1$. In dimension one, Jotsaroop and Shrivastava \cite{JS22} extended these results to the non-Banach range. Maximal and square-function estimates for bilinear Bochner--Riesz operators have been studied in \cite{JL20, JS22, CJSS23, CS23}, and the theory has been extended to multilinear Bochner--Riesz operators in \cite{Shuin24, HLZ24}.
	
	At the critical index $\alpha_c=n-\frac12$, strong-type $L^{p_1}\times L^{p_2}\to L^p$ estimates are known for all $1\leq p_1,p_2\leq\infty$ except $(p_1,p_2)=(1,1)$, see \cite{BGSY15,JSS21}. At this endpoint, Bernicot et al. \cite[Proposition 4.2]{BGSY15} showed that the strong-type estimate fails. Jotsaroop, Shrivastava and Shuin \cite{JSS21} proved weighted estimates at the critical index for bilinear Muckenhoupt weights when $1<p_1,p_2\leq\infty$. More recently, Carro, Luque and S\'anchez-Pascuala \cite{CLS26} obtained weighted endpoint estimates for Orlicz-type spaces near $L^1$. The question of whether $\mathcal B^{n-\frac12}$ satisfies the weak-type $(1,1;\frac12)$ estimate, the bilinear analogue of Christ's theorem \cite{Christ88}, has remained open for a long time.
    
    \subsection{Endpoint estimates for bilinear operators} For bilinear operators, endpoint weak-type estimates are considerably more delicate than in the linear case. Multilinear Calder\'on--Zygmund theory goes back to Coifman and Meyer \cite{CM75}. Grafakos and Torres \cite{GT02} proved that a multilinear Calder\'on--Zygmund operator whose kernel satisfies the standard size and H\"older-type regularity conditions maps $L^1\times L^1\times\dots L^1$ into $L^{\frac1m,\infty}$. The regularity assumption on the kernel is essential, and the endpoint estimate can fail when it is dropped. 
    
    A striking contrast is the class of bilinear rough singular integrals
	\[T_\Omega(f,g)(x)=\mathrm{p.v.}\int_{\R^n}\int_{\R^n}\frac{\Omega\big((y,z)/|(y,z)|\big)}{|(y,z)|^{2n}}\,f(x-y)\,g(x-z)\,dy\,dz,\]
	where $\Omega$ has mean zero on $\mathbb S^{2n-1}$. Grafakos, He and Honz\'ik \cite{GHH18} proved that $T_\Omega$ maps $L^2\times L^2$ into $L^1$ when $\Omega\in L^\infty(\mathbb S^{2n-1})$. This was extended to $\Omega\in L^q$ with $q>\frac43$ by Grafakos, He and Slav\'ikov\'a \cite{GHS20}, and to the full range $q>1$ and all $L^{p_1}\times L^{p_2}\to L^p$ estimates with $1<p_1,p_2<\infty$ by Dosidis and Slav\'ikov\'a \cite{DS24}. We refer to Bhojak and Shrivastava~\cite{bhojak26} for an alternative proof of these results. At the endpoint, however, the analogy with the linear theory breaks down. Honz\'ik \cite{Honzik26} recently showed that, for every $n\geq1$, there exist bounded mean-zero functions $\Omega$ on $\mathbb S^{2n-1}$ for which $T_\Omega$ fails to map $L^1\times L^1$ into $L^{\frac12,\infty}$. 
	
		Similar contrasts between linear and bilinear endpoint estimates have appeared in other contexts as well. In the linear case, the lacunary spherical maximal function is bounded on $L^p(\R^n)$ for all $p>1$ and $n\geq2$ \cite{Calderon79}, and Christ \cite{Christ88} proved that it maps $H^1(\R^n)$ into $L^{1,\infty}(\R^n)$. In contrast, Lee, Li, Wen and Wu \cite{LLWW26} recently showed that the bilinear lacunary spherical maximal function fails to map $H^1(\R^n)\times H^1(\R^n)$ into $L^{\frac{1}{2},\infty}(\R^n)$ for every $n\geq1$. Similarly, the bilinear Hilbert transform is bounded in a large range of exponents by the work of Lacey and Thiele \cite{LT97, LT99}, but Alves and Grafakos \cite{AG26} showed that it fails to satisfy even the restricted weak-type $(1,1;\frac12)$ estimate.

	These results show that a weak-type endpoint estimate for a linear operator need not have a bilinear counterpart. In particular, the counterexample of Honz\'ik \cite{Honzik26} is relevant here as it is closely related to our problem. Note that at the critical index, the bilinear Bochner--Riesz kernel has exactly the same decay $|(y,z)|^{-2n}$ as the kernel of $T_{\Omega}$. The main idea behind Honz\'ik's counterexample is to choose $\Omega$ so that the contributions from different spatial scales have the same sign. On the other hand, the leading term of the bilinear Bochner--Riesz kernel oscillates in the radial variable. Indeed, at the critical index, the asymptotic expansion of the Bessel function gives
	\begin{equation}\label{kernelasy}
		K^{n-\frac12}(y,z)=\frac{c_+e^{2\pi i|(y,z)|}+c_-e^{-2\pi i|(y,z)|}}{|(y,z)|^{2n}}+O\bigl(|(y,z)|^{-2n-1}\bigr)
	\end{equation}
	as $|(y,z)|\to\infty$, where $c_+$ and $c_-$ are constants. Although the main term has the same size as a rough singular integral kernel, the factor $e^{\pm2\pi i|(y,z)|}$ prevents contributions from different scales from piling up and provides additional cancellation. Our main result shows that this oscillation is enough to restore the endpoint estimate.
    %%%%%%%%%%%%%%%%%%%%%%
	\subsection{Main result}
%%%%%%%%%%%%%%%%%%%%%%%%%%
The following theorem is the main result of the paper. 
	\begin{theorem}\label{thm:main}
		Let $n\geq1$. Then the bilinear Bochner--Riesz operator $\mathcal B^{n-\frac12}$ is bounded from $L^1(\R^n)\times L^1(\R^n)$ to $L^{\frac12,\infty}(\R^n)$. In particular,
		\[\|\mathcal B^{n-\frac{1}{2}}(f,g)\|_{L^{\frac{1}{2},\infty}(\R^n)}\leq C_n\|f\|_{L^1(\R^n)}\|g\|_{L^1(\R^n)}.\]
	\end{theorem}
	We make several comments on Theorem \ref{thm:main}.
	\begin{enumerate}
		\item Theorem \ref{thm:main} is the bilinear analogue of Christ's weak-type $(1,1)$ estimate for $B^{\frac{n-1}{2}}$ \cite{Christ88}. In view of \cite[Proposition 4.2]{BGSY15}, it is sharp in the sense that $L^{\frac12,\infty}$ cannot be replaced by $L^{\frac12}$. Together with the strong-type estimates away from $(1,1)$, it gives a complete description of the behaviour of $\mathcal B^{n-\frac12}$ at $(1,1)$.
 
		\item In view of Honz\'ik's counterexample \cite{Honzik26}, Theorem \ref{thm:main} shows that, for bilinear kernels of critical size $|(y,z)|^{-2n}$, the endpoint weak-type estimate is decided by the finer structure of the kernel and not by its size. In the bilinear Bochner--Riesz case, the radial oscillation in \eqref{kernelasy} supplies the required cancellation.
        \item The proof of Theorem~\ref{thm:main} require new and noble ideas.
		\item We highlight the new features of the proof. To the best of our knowledge, Theorem \ref{thm:main} is the first positive $L^1\times L^1\to L^{\frac12,\infty}$ estimate for a bilinear operator whose kernel has the critical size $|(y,z)|^{-2n}$, but fails the standard regularity conditions of bilinear Calder\'on--Zygmund theory. The recent results for rough or singular bilinear operators at this endpoint have been negative \cite{Honzik26, LLWW26, AG26}. 
        \item The main technical novelty is the bilinear $TT^*$ argument. Christ's linear argument relies on the fact that the correlation kernels are convolution kernels. In the bilinear setting this structure is lost: the correlation kernels $\mathfrak K_{j,k}$ defined in \eqref{eq:Corelation kernel} are integrals of kernels on $\R^{2n}$ over the $n$-dimensional diagonal, and their oscillatory behaviour must be analysed directly. The estimates for $\mathfrak K_{j,k}$ and the splitting of the bad--bad term into near-diagonal and off-diagonal parts are new. The method works uniformly in all dimensions $n\geq1$, using angular stationary phase when $n\geq2$ and van der Corput's lemma when $n=1$. We expect it to apply to endpoint problems for other bilinear operators with oscillatory kernels.
	\end{enumerate}
	\subsection{Strategy of the proof}
	Our proof is motivated by Christ's method for the linear Bochner--Riesz operator at the critical index \cite{Christ88}. The main idea is to combine the Calder\'on--Zygmund decomposition with estimates that exploit the oscillation of the kernel. We now briefly	describe the main steps.	
	\begin{enumerate}
		\item We first apply the Calder\'on--Zygmund decomposition to both $f$ and $g$, writing each as a good function in $L^1\cap L^\infty$ plus a sum of bad functions with mean zero, supported on pairwise disjoint dyadic cubes. The three terms containing at least one good function are controlled using the mixed $L^1\times L^\infty\to L^1$ estimates. The main difficulty is therefore the term containing two bad functions. 
		\item To estimate this term, we decompose the Bochner--Riesz kernel into pieces $K_j$ supported in dyadic annuli $|(y,z)|\sim2^j$ in $\R^{2n}$. The corresponding operators are
		\[T_j(f,g)(x)=\int_{\R^n}\int_{\R^n}K_j(x-y,x-z)\,f(y)\,g(z)\,dy\,dz.\]
		As in \cite{Christ88}, for each $j$ we group the bad cubes according to the ratio of their side length to the kernel scale $2^j$. The bad cubes that are large compared with $2^j$ are removed by discarding a suitable enlargement of their union. 
		
		\item We next split the bad--bad term into two parts, according to whether the variables $y$ and $z$ are close to each other or well separated relative to the scale $2^j$. For the part close to the diagonal $y=z$, we obtain an $L^1$ estimate (Proposition \ref{prop:mainestimateL1}) using only the support properties of the kernel and the disjointness of the bad cubes.
	
		\item For the part away from the diagonal, we use the cancellation of the bad functions together with a bilinear $TT^*$ argument. Expanding the square of the $L^2$ norm of $\sum_j T_j(b_1,b_2)$ leads to correlation kernels of the form
        \begin{equation}\label{eq:Corelation kernel}
            \mathfrak K_{j,k}(y,z,y',z')=\int_{\R^n}K_j(x-y,x-z)\,\overline{K_k(x-y',x-z')}\,dx.
        \end{equation}
		In the linear case, the analogous kernel is $\int_{\R^n}K_j(x-y)\overline{K_k(x-y')}\,dx$, a convolution kernel on $\R^n$. It can be handled by Plancherel's theorem and the explicit behaviour of the multiplier near the unit sphere. In the bilinear case, $K_j$ and $K_k$ are functions on $\R^{2n}$, but the integration takes place only along the $n$-dimensional diagonal $\{(x,x):x\in\R^n\}$. Consequently, $\mathfrak K_{j,k}$ is not a convolution kernel on $\R^{2n}$, and the corresponding linear estimates cannot be applied directly.
		
		\item To overcome this difficulty, we prove oscillatory estimates for these kernels $\mathfrak K_{j,k}$ directly; see Section \ref{sec:kernel}. By \eqref{kernelasy}, $\mathfrak K_{j,k}$ is essentially an oscillatory integral over $\R^n$ whose phase is a difference of two distance functions in $\R^{2n}$, restricted to the diagonal. When $n\geq2$, we use polar coordinates and angular stationary phase, following the ideas of Christ \cite{Christ88}. When $n=1$, we use the first- and second-derivative forms of van der Corput's lemma. Combined with the cancellation of the bad functions, these estimates give an $L^2$ bound (Proposition \ref{prop:mainestimateL2}) with geometric decay in the relative scale between the bad cubes and the kernel.
		
		\item Finally, we combine the $L^1$ estimate for the part near the diagonal with the $L^2$ estimate for the remaining part, using Chebyshev's inequality at the appropriate heights. The geometric decay in the $L^2$ estimate allows us to sum over all relative scales, which completes the proof of Theorem \ref{thm:main}.
	\end{enumerate}
	
	\subsection*{Organization of the paper}
	In Section \ref{sec:CZ}, we apply the Calder\'on--Zygmund decomposition and reduce the proof of Theorem \ref{thm:main} to estimating the bad--bad term. In Section \ref{sec:bad}, we decompose this term into a part near the diagonal and a part away from it, and reduce the desired estimate to an $L^1$ estimate (Proposition \ref{prop:mainestimateL1}) and an $L^2$ estimate (Proposition \ref{prop:mainestimateL2}). We prove these two estimates in Sections \ref{sec:L1} and \ref{sec:L2}, respectively. Finally, in Section \ref{sec:kernel}, we prove the oscillatory estimates for the correlation kernels used in the $L^2$ argument.
	
	\section{Calder\'on--Zygmund decomposition}\label{sec:CZ}
	By homogeneity, it is enough to consider $\|f\|_1=\|g\|_1=1$. Let $\lambda>0$. We apply the dyadic Calder\'on--Zygmund decomposition to $f$ and $g$ at height $\lambda^{\frac{1}{2}}$ to get
	\[f=f_0+b,\quad b=\sum_Qb_Q,\qquad \text{and}\qquad g=g_0+c,\quad c=\sum_Pc_P,\]
	where the cubes $Q$ (respectively $P$) are pairwise disjoint and satisfy the following properties.
	\begin{align*}
		\sum_Q|Q|+\sum_P|P|&\lesssim\lambda^{-\frac{1}{2}},\\
		\operatorname{supp}b_Q\subset Q,\quad \int b_Q=0,\quad\|b_Q\|_1&\lesssim\lambda^{\frac{1}{2}}|Q|,\\
		\operatorname{supp}c_P\subset P,\quad \int c_P=0,\quad\|c_P\|_1&\lesssim\lambda^{\frac{1}{2}}|P|,\\
		\|f_0\|_\infty+\|g_0\|_\infty&\lesssim\lambda^{\frac{1}{2}},\\
		\|f_0\|_2^2+\|g_0\|_2^2&\lesssim\lambda^{\frac{1}{2}}.
	\end{align*}
	We can organize the bad functions based on the length of the cubes. For $k\in\Z$, define
	\[b_k=\sum_{\ell(Q)=2^k}b_Q,\qquad c_k=\sum_{\ell(P)=2^k}c_P,\]
	where $\ell(Q)$ is the side length of $Q$. Consider the set
	\[E=\bigcup_Q5Q\cup\bigcup_P5P.\]
	Then, from the above conditions, we get
	\begin{equation}\label{exceptionalset}
		|E|\lesssim\lambda^{-\frac{1}{2}}.
	\end{equation}
	By Calder\'on--Zygmund decomposition, we have
	\begin{align*}
		&|\{x:|\mathcal B^{n-\frac{1}{2}}(f,g)(x)|>\lambda\}|\\
		\leq& |E|+|\{x\notin E:|\mathcal B^{n-\frac{1}{2}}(f_0,g_0)(x)|>\frac{\lambda}{4}\}|+|\{x\notin E:|\mathcal B^{n-\frac{1}{2}}(b,g_0)(x)|>\frac{\lambda}{4}\}|\\
		&+|\{x\notin E:|\mathcal B^{n-\frac{1}{2}}(f_0,c)(x)|>\frac{\lambda}{4}\}|+|\{x\notin E:|\mathcal B^{n-\frac{1}{2}}(b,c)(x)|>\frac{\lambda}{4}\}|.
	\end{align*}
	We first note that $\mathcal B^{n-\frac{1}{2}}$ is bounded from $L^1(\R^n)\times L^\infty(\R^n)$ into $L^1(\R^n)$ for every $n\geq1$. This follows from \cite[Theorem 4.1]{BGSY15} when $n=1$ and \cite[Theorem 4.9]{BGSY15} when $n\geq2$. Using this estimate and Chebyshev's inequality, we control all the terms containing at least one good function in the Calder\'on--Zygmund decomposition.
	\subsection*{Good--good term}
	By $L^1\times L^\infty\to L^1$-bound for $\mathcal B^{n-\frac{1}{2}}$ and Chebyshev's inequality, we get
	\begin{align}\label{gg}
		|\{x:|\mathcal B^{n-\frac{1}{2}}(f_0,g_0)(x)|>\frac{\lambda}{4}\}|\lesssim \lambda^{-1}\|f_0\|_1\|g_0\|_\infty \lesssim\lambda^{-\frac{1}{2}}.
	\end{align}
	\subsection*{Good--bad terms}
	Again by $L^1\times L^\infty\to L^1$-bound for $\mathcal B^{n-\frac{1}{2}}$ and Chebyshev's inequality, we get
	\begin{align}\label{bg}
		|\{x:|\mathcal B^{n-\frac{1}{2}}(b,g_0)(x)|>\frac{\lambda}{4}\}|\lesssim \lambda^{-1}\|b\|_1\|g_0\|_\infty\lesssim \lambda^{-\frac{1}{2}}.
	\end{align}
	
	Similarly, the $L^\infty\times L^1\to L^1$-bound for $\mathcal B^{n-\frac{1}{2}}$ together with Chebyshev's inequality gives
	\begin{align}\label{gb}
		|\{x:|\mathcal B^{n-\frac{1}{2}}(f_0,c)(x)|>\frac{\lambda}{4}\}|\lesssim \lambda^{-1}\|f_0\|_\infty\|c\|_1\lesssim \lambda^{-\frac{1}{2}}.
	\end{align}
	
	\subsection*{Bad--bad term}
	Recall that the kernel of $\mathcal B^{n-\frac{1}{2}}$ is given by
	\[K^{n-\frac{1}{2}}(y,z)=\frac{\Gamma(n+\frac{1}{2})}{\pi^{n-\frac{1}{2}}}\frac{J_{2n-\frac{1}{2}}(2\pi|(y,z)|)}{|(y,z)|^{2n-\frac{1}{2}}}.\]
	From the derivative properties of Bessel functions, we have for every positive integer $N$
	\begin{equation}\label{kernelderi}
		|\partial^\gamma K(y,z)|\lesssim |(y,z)|^{-2n},\qquad |(y,z)|\geq1,\quad |\gamma|\leq N.
	\end{equation}
	
	Let $\psi_0\in C_c^\infty(\R^{2n})$ be a radial function such that $0\leq\psi_0\leq1$, $\psi_0(y,z)=1$ when $|(y,z)|\leq1$ and $\operatorname{supp}\psi_0\subset\{|(y,z)|\leq2\}$. Define $\psi(y,z)=\psi_0(y,z)-\psi_0(2(y,z))$. Then
	\[\psi_0(y,z)+\sum_{j\geq1}\psi(2^{-j}(y,z))=1.\]
	Write $K_0=\psi_0K$, $K_j(y,z)=\psi(2^{-j}(y,z))K(y,z)$ for $j\geq1$. Let $T_j$ denote the bilinear operator associated with the kernel $K_j$,
	\[T_j(f,g)(x)=\int_{\R^n}\int_{\R^n}K_j(x-y,x-z)f(y)g(z)\,dy\,dz.\]
	Using the support of $\psi$, \eqref{kernelasy} and \eqref{kernelderi}, we get for every $j\geq1$
	\begin{align}
		\operatorname{supp}K_j&\subset\{2^{j-1}\leq|X|\leq2^{j+1}\},\nonumber\\
		\|\partial^\gamma K_j\|_\infty&\lesssim_\gamma2^{-2nj},\qquad |\gamma|\leq N,\nonumber\\
		\sup_{y,z\in\R^n}\int_{\R^n}|\partial^\gamma K_j(x-y,x-z)|\,dx&\lesssim_\gamma2^{-nj},\qquad |\gamma|\leq N.\label{deri Kj}
	\end{align}
	
	Since $K_0$ is bounded, it is easy to see that the operator $T_0$ maps $L^{p_1}\times L^{p_2}$ to $L^p$ for all $1\leq p_1,p_2\leq\infty$ satisfying $\frac{1}{p}=\frac{1}{p_1}+\frac{1}{p_2}$. We claim the following estimate holds for $\sum_{j\geq1}T_j$.
	\begin{equation}\label{Tj bb}
		|\{x\notin E:|\sum_{j\geq1}T_j(b,c)(x)|>\frac{\lambda}{8}\}|\lesssim \lambda^{-\frac{1}{2}}.
	\end{equation}
	If we assume \eqref{Tj bb}, then together with the estimates of $T_0$, we have
	\begin{equation}\label{bb}
		|\{x\notin E:|\mathcal B^{n-\frac{1}{2}}(b,c)(x)|>\frac{\lambda}{4}\}|\lesssim \lambda^{-\frac{1}{2}}.
	\end{equation}
	Combining the estimates in \eqref{exceptionalset}, \eqref{gg}, \eqref{bg}, \eqref{gb} and \eqref{bb} proves the Theorem \ref{thm:main}.
	
	Therefore, we need to prove the claim \eqref{Tj bb}.

	\section{Decomposition of the bad--bad term}\label{sec:bad}
	We now prove \eqref{Tj bb} for the bad--bad term. We first organize the bad functions according to the sizes of their supporting cubes relative to the scale of $K_j$. We then divide the interactions into two parts, depending on whether the variables $y$ and $z$ are close to each other or well separated.

	Fix $j\geq1$ and $x\notin E$. Suppose that $y\in Q$, $z\in P$ and $K_j(x-y,x-z)\ne0$. Since $x\notin5Q\cup5P$, we have
	\[|x-y|>2\ell(Q),\qquad |x-z|>2\ell(P),\]
	whereas $|x-y|,|x-z|\leq2^{j+1}$. Consequently, we obtain $\ell(Q)<2^j$ and $\ell(P)<2^j$. Thus we only need to consider cubes of the length $\ell(Q)=2^{j-s}$ and $\ell(P)=2^{j-t}$ for $s,t\geq0$. Using this observation, we can rewrite the operator $\sum\limits_{j\geq1}T_j$ as
	\begin{equation}\label{Tj reorganized}
		\sum_{j\geq1}T_j(b,c)(x)=\sum_{s,t\geq0}H_{s,t}(b,c)(x),
	\end{equation}
	where
	\[H_{s,t}(b,c)(x)=\sum_{j\geq1}T_j(b_{j-s},c_{j-t})(x).\]
    Note that the function $H_{s,t}$ is symmetric in the sense that $H_{s,t}(b,c)=H_{t,s}(c,b)$. Thus, it is enough to estimate $\sum_{s\geq t\geq0}H_{s,t}(b,c)$.
	
	Further, observe that $K_j(x-y,x-z)\ne0$ implies $|y-z|\leq2^{j+2}$. Therefore, for fixed $s,t,j$, it is enough to consider pairs of bad cubes satisfying $\operatorname{dist}(Q,P)\leq 2^{j+2}$. Define
	\begin{equation}\label{Fj}
		F_j(y,z)=\sum_{\substack{\ell(Q)=2^{j-s},\ \ell(P)=2^{j-t}\\\operatorname{dist}(Q,P)\leq 2^{j+2}}}b_Q(y)c_P(z).
	\end{equation}
	Then, for $x\notin E$, we have 
	\begin{align}\label{TjinFj}
		T_j(b_{j-s},c_{j-t})(x)&=\iint K_j(x-y,x-z)F_j(y,z)\,dy\,dz.
	\end{align}
	
	Since each pair $(Q,P)$ in $F_j(y,z)$ occurs once and the rectangles $Q\times P$ are pairwise disjoint in $\R^{2n}$. Thus for $y\in Q$ and $z\in P$, we have  
	\[|y-z|\leq\operatorname{dist}(Q,P)+\sqrt n\ell(Q)+\sqrt n\ell(P)\leq 2^{j+2}+\sqrt n2^{j-s}+\sqrt n 2^{j-t}\leq (2n+4) 2^j.\]
	Consequently, we obtain 
	\begin{align}\label{estimateFj}
		|F_j(y,z)|\leq|b_{j-s}(y)||c_{j-t}(z)|\mathbf 1_{\{|y-z|\leq (2n+4) 2^j\}}.
	\end{align}
	
	Let $\varphi\in C_c^\infty(\R^n)$ such that $0\leq\varphi\leq1$, $\varphi=1$ on $B(0,1)$ and $\operatorname{supp}\varphi\subset B(0,2)$. For every $Q$ with $\ell(Q)=2^{j-s}$, let $h_j=2^{j-\frac{s}{4}}$ and define
	\[\varphi_{Q,j}(z)=\varphi\left(\frac{z-\rho_Q}{(n+1)h_j}\right),\]
	where $\rho_Q$ denotes the center of $Q$. Note that $\varphi_{Q,j}(z)\neq 0$ implies $|z-\rho_Q|\leq 2(n+1)h_j$ whereas  $\varphi_{Q,j}(z)\neq1$ implies $|z-\rho_Q|\geq(n+1)h_j$.
	
	We decompose $F_j=F_j^1+F_j^2$, where
	\begin{align*}
		F_j^1(y,z)&=\sum_{(Q,P)}b_Q(y)c_P(z)\varphi_{Q,j}(z),\\
		F_j^2(y,z)&=\sum_{(Q,P)}b_Q(y)c_P(z)(1-\varphi_{Q,j}(z)).
	\end{align*}
	Consequently, we get that 
	\[H_{s,t}(b,c)=H^{1}_{s,t}(b,c)+H^{2}_{s,t}(b,c).\]
	Since $0\leq \varphi\leq 1$, the estimate in \eqref{estimateFj} continues to hold for both the pieces $F^{1}_j$ and $F^{2}_j$.
	
	We prove the following estimates for $H^{1}_{s,t}(b,c)$ and $H^{2}_{s,t}(b,c)$.
	
	\begin{proposition}\label{prop:mainestimateL1}
		For $s\geq t\geq0$, we have
		\[\|H^1_{s,t}(b,c)\|_{L^1(E^c)}\lesssim2^{-\frac{ns}{4}}\lambda^{\frac{1}{2}}.\]
	\end{proposition}
	
	\begin{proposition}\label{prop:mainestimateL2}
		For $s\geq t\geq0$, we have
		\[\|H_{s,t}^2(b,c)\|_2\lesssim2^{-\frac{s}{8}}\lambda^{\frac{3}{4}}.\]
	\end{proposition}
	We assume both the results above for the moment. Then
	\begin{align*}
		&\left|\left\{x\notin E:|\sum_{j\ge1}T_j(b,c)(x)|>\frac{\lambda}{8}\right\}\right|\\
		&= \left|\left\{x\notin E:|\sum_{s,t\geq 0}H_{s,t}(b,c)(x)|>\frac{\lambda}{8}\right\}\right|\\
		&\lesssim \left|\left\{x\notin E:|\sum_{s\geq t\geq 0}H_{s,t}(b,c)(x)|>\frac{\lambda}{16}\right\}\right|+\left|\left\{x\notin E:|\sum_{s\geq t\geq 0}H_{s,t}(c,b)(x)|>\frac{\lambda}{16}\right\}\right|\\
		&\lesssim \left|\left\{x\notin E:|\sum_{s\geq t\geq 0}H^{1}_{s,t}(b,c)(x)|>\frac{\lambda}{32}\right\}\right|+ \left|\left\{x:|\sum_{s\geq t\geq 0}H^{2}_{s,t}(b,c)(x)|>\frac{\lambda}{32}\right\}\right|\\
		&\hspace{.5cm}+\left|\left\{x\notin E:|\sum_{s\geq t\geq 0}H^{1}_{s,t}(c,b)(x)|>\frac{\lambda}{32}\right\}\right|+ \left|\left\{x:|\sum_{s\geq t\geq 0}H^{2}_{s,t}(c,b)(x)|>\frac{\lambda}{32}\right\}\right|\\
		&\lesssim \sum_{s\geq t\geq 0} \left(\left|\left\{x\notin E:|H^{1}_{s,t}(b,c)(x)|>\lambda_{s,t}\right\}\right|+\left|\left\{x:|H^{2}_{s,t}(b,c)(x)|>\lambda_{s,t}\right\}\right|\right)\\
		&\hspace{.5cm}+ \sum_{s\geq t\geq 0}\left(\left|\left\{x\notin E:|H^{1}_{s,t}(c,b)(x)|>\lambda_{s,t}\right\}\right|+\left|\left\{x:|H^{2}_{s,t}(c,b)(x)|>\lambda_{s,t}\right\}\right|\right),
	\end{align*}
	where $\lambda_{s,t}=c_0\lambda2^{-\frac{s}{16}}$ with $c_0>0$ a fixed constant such that  $\sum_{s\geq t\ge0}\lambda_{s,t}\leq \frac{\lambda}{16}$. Applying  Chebyshev's inequality along with estimates from Proposition \ref{prop:mainestimateL1} and Proposition \ref{prop:mainestimateL2}, we get that 
	\begin{align*}
		\left|\left\{x\notin E:|\sum_{j\ge1}T_j(b,c)(x)|>\frac{\lambda}{8}\right\}\right|\lesssim&  \sum_{s\geq t\geq 0}\left(\frac{\|H^{1}_{s,t}(b,c)\|_{L^1(E^c)}}{\lambda_{s,t}}+\frac{\|H^{2}_{s,t}(b,c)\|_2^2}{\lambda_{s,t}^2}\right)\\
        &+\sum_{s\geq t\geq 0}\left(\frac{\|H^{1}_{s,t}(c,b)\|_{L^1(E^c)}}{\lambda_{s,t}}+\frac{\|H^{2}_{s,t}(c,b)\|_2^2}{\lambda_{s,t}^2}\right)\\
		\lesssim& \sum_{s\geq t\geq 0}\left(\frac{2^{-\frac{ns}{4}}\lambda^{\frac{1}{2}}}{2^{-\frac{s}{16}}\lambda}+\frac{2^{-\frac{s}{4}}\lambda^{\frac{3}{2}}}{2^{-\frac{s}{8}}\lambda^2}\right)\\
		\lesssim& \lambda^{-\frac{1}{2}}.
	\end{align*}
	Therefore, the estimate in \eqref{Tj bb} is reduced to proving Propositions \ref{prop:mainestimateL1} and \ref{prop:mainestimateL2}.  We first prove  Proposition \ref{prop:mainestimateL1}. 
	
	\section{The near-diagonal $L^1$ estimate}\label{sec:L1}
	\begin{proof}[Proof of Proposition \ref{prop:mainestimateL1}]
		Note that from \eqref{deri Kj}, we have
		\begin{equation}\label{singlejprop1}
			\left\|\int\!\!\int K_j(\cdot-y,\cdot-z)F_j^1(y,z)\,dy\,dz\right\|_1\lesssim2^{-nj}\|F_j^1\|_1.
		\end{equation}
		Therefore, we need to estimate the $L^1$ norm of $F_j^1$. Using \eqref{estimateFj} and the definition of $F^{1}_j$, we have
		\begin{align*}
			\|F^{1}_j\|_1
			&\leq \sum_{(Q,P)}\int_Q\int_P|b_Q(y)||c_P(z)|\varphi_{Q,j}(z)\,dy\,dz\\
			&\leq \sum_{(Q,P)}\|b_Q\|_1\int|c_P(z)|\varphi_{Q,j}(z) dz\\
			&\leq \int|c_{j-t}(z)|\left(\sum_Q\|b_Q\|_1\varphi_{Q,j}(z)\right)dz .
		\end{align*}
		
		For a fixed $z\in\R^n$, consider the collection of intervals $Q$ given by 
		\[\mathcal Q(z)=\{Q:\ \ell(Q)=2^{j-s},\ \varphi_{Q,j}(z)\neq 0\}.\]
		We can see that $\varphi_{Q,j}(z)\neq0$ implies $|z-\rho_Q|\leq 2(n+1)h_j$, therefore every $Q\in\mathcal Q(z)$ has its center in $B(z,Ch_j)$. Since cubes in $\mathcal Q(z)$ are pairwise disjoint and $\ell(Q)=2^{j-s}$, we have
		\[\sum_{Q\in\mathcal Q(z)}|Q|  \lesssim h_j^n+2^{n(j-s)}\lesssim h_j^n.\]
		Therefore, we get that 
		\begin{align*}
			\|F^{1}_j\|_1
			&\leq \int|c_{j-t}(z)|\left(\sum_Q\|b_Q\|_1\varphi_{Q,j}(z)\right)dz \\
			&\lesssim  \int|c_{j-t}(z)|\left(\sum_{Q\in\mathcal Q(z)}\|b_Q\|_1\right)\,dz\\
			&\lesssim \lambda^{\frac{1}{2}}\int|c_{j-t}(z)|\left(\sum_{Q\in\mathcal Q(z)}\|b_Q\|_1\right)\,dz\\
			& \lesssim \lambda^{\frac{1}{2}}h_j^n\|c_{j-t}\|_1.
		\end{align*}
		Substituting the above estimate in \eqref{singlejprop1}, we get
		\[\left\|\iint K_j(\cdot-y,\cdot-z)F^{1}_j(y,z)\,dy\,dz\right\|_1\lesssim \lambda^{\frac{1}{2}} \frac{h_j^n}{2^{nj}}\|c_{j-t}\|_1 \lesssim 2^{-\frac{ns}{4}}\lambda^{\frac{1}{2}}\|c_{j-t}\|_1 .\]
		Now, summing in $j$, we obtain
		\[\|H^{1}_{s,t}(b,c)\|_{L^1(E^c)}\lesssim2^{-\frac{ns}{4}}\lambda^{\frac{1}{2}}.\]
		This completes the proof of Proposition \ref{prop:mainestimateL1}.
	\end{proof}
	
	\section{The $L^2$ estimate for the separated term}\label{sec:L2}
	To estimate $H^2_{s,t}(b,c)$, we follow ideas from Christ \cite{Christ88} and use a $TT^*$ argument that combines the cancellation of the bad functions with the oscillation of the kernel. We first obtain some local $L^1$ estimates for the bad functions.
	\begin{lemma}\label{lem:badterm}
		Let $B\subset\R^n$ be a ball of radius $r>0$. Then for every $k$,
		\begin{equation}\label{badtermestimate}
			\int_B|b_k(y)|\,dy\lesssim\lambda^{\frac{1}{2}}(r+2^k)^n,\qquad
			\int_B|c_k(z)|\,dz\lesssim\lambda^{\frac{1}{2}}(r+2^k)^n.
		\end{equation}
		Moreover,
		\begin{equation}\label{estisumFj}
			\sum_{j\geq1}2^{-nj}\|F_j\|_1\lesssim\lambda^{\frac{1}{2}}.
		\end{equation}
	\end{lemma}
	\begin{proof}
		If a cube $Q$ of side length $2^k$ intersects $B$, then it is contained in the $C2^k$-neighborhood of $B$. Since the cubes $Q$ are disjoint and are of sidelength $2^k$, we have
		\begin{align*}
			\int_B|b_k(y)|\,dy&\leq \sum_{Q\cap B\ne\emptyset}\int |b_Q(y)|\,dy\\
			&\lesssim \lambda^{\frac{1}{2}}\sum_{Q\cap B\ne\emptyset}|Q|\\
			&\lesssim \lambda^{\frac{1}{2}}(r+2^k)^n.
		\end{align*}
		Similarly, we can obtain the estimate for $c_k$. By \eqref{estimateFj} and \eqref{badtermestimate} with a ball of radius $(2n+4)2^j$, we get
		\begin{align*}
			\|F_j\|_1&\leq\int|b_{j-s}(y)|\int_{|y-z|\leq (2n+4) 2^j}|c_{j-t}(z)|\,dz\,dy\\
			&\lesssim\lambda^{\frac{1}{2}}2^{nj}\|b_{j-s}\|_1.
		\end{align*}
		Thus
		\begin{align*}
			\sum_{j\geq1}2^{-nj}\|F_j\|_1\lesssim \sum_{j\geq1}\lambda^{\frac{1}{2}}\|b_{j-s}\|_1\leq \lambda^{\frac{1}{2}}.
		\end{align*}
	\end{proof}
	
	\begin{proof}[Proof of Proposition \ref{prop:mainestimateL2}]
		We exploit the cancellation of the bad functions to obtain a suitable representation of $H_{s,t}^2(b,c)$.
		\subsection*{Cancellation}
		On the support of $F_j^2$, we have $1-\varphi_{Q,j}(z)\neq0$, and hence
		$|z-\rho_Q|\geq(n+1)h_j$. Since $\ell(Q)=2^{j-s}$, every $y\in Q$ satisfies
		\[|y-\rho_Q|\leq\frac{\sqrt n}{2}2^{j-s}\leq\frac{\sqrt n}{2}h_j.\]
		Consequently, we get
		\begin{equation}\label{diff away}
			|y-z|\geq|z-\rho_Q|-|y-\rho_Q|\geq\frac32h_j.
		\end{equation}
		
		We now consider the case $j<s$. Fix an admissible pair $(Q,P)$ and $z$ such that $1-\varphi_{Q,j}(z)\neq0$. Since $\int b_Q=0$, the fundamental theorem of calculus gives
		\begin{align*}
			\int b_Q(y)K_j(x-y,x-z)\,dy
			&=\int b_Q(y)\bigl(K_j(x-y,x-z)-K_j(x-\rho_Q,x-z)\bigr)\,dy\\
			&=-\int_0^1\int b_Q(y)(y-\rho_Q)\cdot\nabla_1K_j(x-\rho_Q-\vartheta(y-\rho_Q),x-z)\,dy\,d\vartheta,
		\end{align*}
		where $\nabla_1K_j$ denotes the gradient in the first $n$ variables. For fixed $0<\vartheta\leq1$, substituting $y=\rho_Q+\frac{w-\rho_Q}{\vartheta}$, we obtain
		\begin{equation}\label{cancellation}
			\int b_Q(y)K_j(x-y,x-z)\,dy=2^{j-s}\int_0^1\int b_Q^\vartheta(w)\cdot\nabla_1K_j(x-w,x-z)\,dw\,d\vartheta,
		\end{equation}
		where the vector-valued function $b_Q^\vartheta$ is defined by
		\[b_Q^\vartheta(w)=-\frac{w-\rho_Q}{2^{j-s}\vartheta^{n+1}}b_Q\left(\rho_Q+\frac{w-\rho_Q}{\vartheta}\right)\mathbf1_{\rho_Q+\vartheta(Q-\rho_Q)}(w).\]
		From the definition of $b_Q^\vartheta$, we have
		\[\operatorname{supp}b_Q^\vartheta\subseteq\rho_Q+\vartheta(Q-\rho_Q)\subseteq Q.\]
		Changing variables from $w$ back to $y$, we also have
		\[\|b_Q^\vartheta\|_1= 2^{-(j-s)}\int_Q|y-\rho_Q|\,|b_Q(y)|\,dy\lesssim \|b_Q\|_1\lesssim \lambda^{\frac12}|Q|,\]
		uniformly in $\vartheta\in(0,1]$. Here, $|\cdot|$ of a vector-valued function denotes its Euclidean norm. 
		
		We define the analogous notation for the modified functions. Let $b_{j-s}^\vartheta=\sum_Q b_Q^\vartheta$, where the sum is over cubes $Q$ of side length $2^{j-s}$. Note that $b_{j-s}^\vartheta$ is vector-valued. Define
		\[F_j^{2,\vartheta}(y,z)=
		\begin{cases}
			\sum_{(Q,P)}b_Q^\vartheta(y)c_P(z)(1-\varphi_{Q,j}(z)),&j<s,\\
			F_j^2(y,z),&j\geq s.
		\end{cases}\]
		Thus, $F_j^{2,\vartheta}$ is vector-valued when $j<s$ and
		scalar-valued when $j\geq s$. Since
		$\operatorname{supp}b_Q^\vartheta\subseteq Q$ and $\|b_Q^\vartheta\|_1\lesssim\|b_Q\|_1$ uniformly in $\vartheta$. Thus the estimates \eqref{estimateFj}, \eqref{badtermestimate}, and \eqref{estisumFj} hold for the modified functions, with constants independent of $\vartheta$. In particular, summing over the admissible pairs, we obtain
		\[\|F_j^{2,\vartheta}\|_1\lesssim\|F_j^2\|_1\leq\|F_j\|_1.\]
		
		We now combine the cases $j<s$ and $j\geq s$. Recall that
		\[H_{s,t}^2(b,c)(x)=\sum_{j\geq1}\iint K_j(x-y,x-z)F_j^2(y,z)\,dy\,dz.\]
		When $j<s$, we use the cancellation identity \eqref{cancellation} to replace $K_j$ by its gradient in the first variable. When $j\geq s$, we keep the original kernel $K_j$. To write both cases in a single expression, let $\gamma_j=\min\{1,2^{j-s}\}$. For each $0<\vartheta\leq1$, define
		\[V_j^\vartheta(x)=
		\begin{cases}
			\gamma_j\iint F_j^{2,\vartheta}(y,z)\cdot\nabla_1K_j(x-y,x-z)\,dy\,dz,&j<s,\\
			\gamma_j\iint K_j(x-y,x-z)F_j^{2,\vartheta}(y,z)\,dy\,dz,&j\geq s.
		\end{cases}\]
		Then, we can write
		\[H_{s,t}^2(b,c)(x)=\sum_{j\geq1}\int_0^1V_j^\vartheta(x)\,d\vartheta.\]
		
		From \eqref{deri Kj} and the estimates for $F_j^{2,\vartheta}$, we have
		\begin{align*}
			\int_{\R^n}\int_0^1\sum_{j\geq1}|V_j^\vartheta(x)|\,d\vartheta\,dx&\lesssim \sum_{j\geq1}\gamma_j2^{-nj}\int_0^1\|F_j^{2,\vartheta}\|_1\,d\vartheta\\
			&\lesssim \sum_{j\geq1}2^{-nj}\|F_j\|_1 \lesssim\lambda^{\frac12},
		\end{align*}
		where the last inequality follows from \eqref{estisumFj}. Therefore, Fubini's
		theorem gives
		\[H_{s,t}^2(b,c)(x)=\int_0^1\sum_{j\geq1}V_j^\vartheta(x)\,d\vartheta\]
		for almost every $x\in\R^n$.
		
		By Minkowski's inequality, we have
		\begin{align*}
			\|H_{s,t}^2(b,c)\|_2&\leq \int_0^1\left\|\sum_{j\geq1}V_j^\vartheta\right\|_2\,d\vartheta\\
			&\leq\sup_{0<\vartheta\leq1}\left\|\sum_{j\geq1}V_j^\vartheta\right\|_2.
		\end{align*}
		Thus, it remains to obtain an $L^2$ estimate for $\sum\limits_{j\geq1}V_j^\vartheta$, uniformly in $\vartheta$.
		\begin{proposition}\label{prop:reduction} The following estimate holds for $V^\vartheta_j$ uniformly in $\vartheta$.
			\begin{align}\label{propvj}
				\left\|\sum_{j\ge1}V^\vartheta_j\right\|_2\lesssim 2^{-\frac{s}{8}}\lambda^{\frac{3}{4}}.
			\end{align}
		\end{proposition}
		If we assume that Proposition \ref{prop:reduction} holds true, then the proof of Proposition \ref{prop:mainestimateL2} is complete.
	\end{proof}
	
	Before giving the proof of Proposition \ref{prop:reduction}, we apply a change of variables to $V_j^\vartheta$.
	\subsection*{Change of variables}
	Let $\Psi:\R^{2n}\to\R^{2n}$ be given by
	\begin{equation}\label{changeofvariable}
		\Psi(y,z)=\left(\frac{y+z}{\sqrt2},\frac{y-z}{\sqrt2}\right)=(u,v),\qquad (y,z)\in\R^n\times\R^n.
	\end{equation}
	Then
	\[\Psi(x-y,x-z)=\left(\frac{2x-y-z}{\sqrt2},\frac{z-y}{\sqrt2}\right)=(\sqrt2x,0)-\Psi(y,z).\]
	Since $K_j$ is radial and $\Psi$ is orthogonal, we have
	\[K_j(x-y,x-z)=K_j\left(\Psi(x-y,\ x-z)\right)=K_j\left((\sqrt2x,0)-\Psi(y,z)\right).\]
	Since the Jacobian of $\Psi$ has absolute value one, applying a change of variable $\Psi(y,z)=(u,v)$, we get
	\[\iint K_j(x-y,x-z)F_j(y,z)\,dy\,dz=\iint K_j\left((\sqrt2x,0)-(u,v)\right)\widetilde F_j(u,v)\,du\,dv,\]
	where $\widetilde{F_j}=F_j\circ\Psi^{-1}$, i.e. $\widetilde{F_j}(u,v)=F_j\left(\frac{u+v}{\sqrt2},\ \frac{u-v}{\sqrt2}\right)$. Moreover, from \eqref{diff away}, we have
	\begin{equation}\label{separation}
		|v|=\frac{|y-z|}{\sqrt2}\geq\frac{3}{2\sqrt2}h_j>h_j.
	\end{equation}
	
	With this change of variables, we first obtain the following estimate for $\widetilde{F_j}$.
	\begin{lemma}\label{lemmalocalesti}
		For every $u_0\in\R^n$, we have
		\begin{align}
			\int_{|u-u_0|\leq r}\int_{\R^n}|\widetilde F_j(u,v)|\,dv\,du&\lesssim \lambda2^{nj}(r+2^{j-s})^n,\qquad 0<r\lesssim 2^j,\label{localmass}\\
			\int_{|u-u_0|\leq r}\int_{\R^n}
			|\widetilde F_j(u,v)|\,dv\,du&\lesssim \lambda2^{nj}r^n,\qquad r\gtrsim 2^j.\label{offdialemmaesti}
		\end{align}
		Moreover, $\widetilde F_j(u,v)=0$ unless $|v|\lesssim2^j$.
	\end{lemma}
	
	\begin{proof}
		The support condition of $F_j$ in \eqref{estimateFj} implies that $|y-z|\leq (2n+4) 2^j$ whenever $F_j(y,z)\neq0$. Since $v=\frac{y-z}{\sqrt2}$, we obtain
		$|v|\lesssim2^j$ on the support of $\widetilde F_j$.
		
		By change of variable $y=\frac{u+v}{\sqrt2}$ and $z=\frac{u-v}{\sqrt2}$, we have
		\[\int_{|u-u_0|\leq r}\int_{\R^n}|\widetilde F_j(u,v)|\,dv\,du
		=\iint_S|F_j(y,z)|\,dy\,dz,\]
		where
		\[S=\bigl\{(y,z):|y+z-\sqrt2u_0|\leq\sqrt2r\bigr\}.\]
		For $(y,z)\in S$ with $F_j(y,z)\neq0$, we have
		\[|2z-\sqrt2u_0|=|(y+z-\sqrt2u_0)-(y-z)|\leq\sqrt2r+(2n+4) 2^j.\]
		
		First consider $r\lesssim2^j$. Then $|2z-\sqrt2u_0|\lesssim 2^j$. Let
		\[B_0=\bigl\{z:|2z-\sqrt2u_0|\lesssim 2^j\bigr\},\qquad B_z=\bigl\{y:|y+z-\sqrt2u_0|\leq\sqrt2r\bigr\}.\]
		For each $z\in B_0$, the corresponding variable $y$ in $S$ lies in the ball $B_z$ with radius $\sqrt2r$. Using \eqref{badtermestimate}, we obtain
		\begin{align*}
			\iint_S|F_j(y,z)|\,dy\,dz&\leq \int_{B_0}|c_{j-t}(z)|\left(\int_{B_z}|b_{j-s}(y)|\,dy\right)dz\\
			&\lesssim\lambda^{\frac12}(r+2^{j-s})^n\int_{B_0}|c_{j-t}(z)|\,dz\\
			&\lesssim\lambda^{\frac12}(r+2^{j-s})^n\lambda^{\frac12}(2^j+2^{j-t})^n\\
			&\lesssim\lambda2^{nj}(r+2^{j-s})^n.
		\end{align*}
		This proves \eqref{localmass}.
		
		Next, when $r\gtrsim 2^j$, for $(y,z)\in S$ with $F_j(y,z)\neq0$, we have
		\begin{align*}
			|2y-\sqrt2u_0|&=|(y+z-\sqrt2u_0)+(y-z)|\\
			&\leq|y+z-\sqrt2u_0|+|y-z|\\
			&\lesssim r+2^j\lesssim r.
		\end{align*}
		Let
		\[B_0=\bigl\{y:|2y-\sqrt2u_0|\lesssim r\bigr\},\qquad B_y=\bigl\{z:|z-y|\leq (2n+4) 2^j\bigr\}.\]
		Then
		\[S\cap\operatorname{supp}F_j\subseteq\bigl\{(y,z):y\in B_0,\ z\in B_y\bigr\}.\]
		Using \eqref{badtermestimate}, we obtain
		\begin{align*}
			\iint_S|F_j(y,z)|\,dy\,dz&\leq \int_{B_0}|b_{j-s}(y)|
			\left(\int_{B_y}|c_{j-t}(z)|\,dz\right)dy\\
			&\lesssim\lambda^{\frac12}(2^j+2^{j-t})^n\int_{B_0}|b_{j-s}(y)|\,dy\\
			&\lesssim\lambda^{\frac12}2^{nj}\lambda^{\frac12}(r+2^{j-s})^n\\
			&\lesssim\lambda2^{nj}r^n.
		\end{align*}
		This proves \eqref{offdialemmaesti}.
	\end{proof}
	
	The estimates in the above lemma also hold for $\widetilde F_j^{2,\vartheta}=F_j^{2,\vartheta}\circ\Psi^{-1}$, uniformly in $\vartheta$, as they have the same support properties and satisfy the same local $L^1$ estimates. In particular, $\widetilde F_j^{2,\vartheta}$ satisfies \eqref{localmass} and \eqref{offdialemmaesti}.
	
	Finally, we rewrite $V_j^\vartheta$ after this change of variables. Since
	\[\nabla_1K_j(y,z)=\frac1{\sqrt2}(\nabla_u+\nabla_v)K_j(\Psi(y,z)),\]
	we define
	\[\mathcal K_j(u,v)=
	\begin{cases}
		\frac1{\sqrt2}(\nabla_u+\nabla_v)K_j(u,v),&j<s,\\
		K_j(u,v),&j\geq s.
	\end{cases}\]
	Thus $\mathcal K_j$ and $\widetilde F_j^{2,\vartheta}$ are both vector-valued when $j<s$ and both scalar-valued when $j\geq s$. Therefore
	\begin{equation}\label{mainreduction}
		V_j^\vartheta(x)=\gamma_j\int_{\R^{2n}}\mathcal K_j\bigl((\sqrt2x,0)-(u,v)\bigr)\cdot
		\widetilde F_j^{2,\vartheta}(u,v)\,du\,dv,
	\end{equation}
	where the `$\cdot$' denotes the Euclidean scalar product in the vector-valued case and ordinary multiplication in the scalar-valued case.
	
	For every $j\geq1$, let
	\[d_j=\begin{cases}n,&j<s,\\1,&j\geq s.\end{cases}\]
	We can write $\mathcal K_j=(\mathcal K_{j,m})_{m=1}^{d_j}$, where $\mathcal K_{j,m}$ is the $m-$th component given by
	\[\mathcal K_{j,m}(u,v)=
	\begin{cases}
		\frac{1}{\sqrt2}(\partial_{u_m}+\partial_{v_m})K_j(u,v),&j<s,\quad 1\leq m\leq n,\\
		K_j(u,v),&j\geq s,\quad m=1.
	\end{cases}\]
	Similarly, we can write $\widetilde F_j^{2,\vartheta}=(\widetilde F_{j,m}^{2,\vartheta})_{m=1}^{d_j}$. Then \eqref{mainreduction} can be written as
	\begin{equation}\label{mainreductioncomponents}
		V_j^\vartheta(x)=\gamma_j\sum_{m=1}^{d_j}\iint_{\R^{2n}}\mathcal K_{j,m}(\sqrt2x-u,-v)\widetilde F_{j,m}^{2,\vartheta}(u,v)\,du\,dv.
	\end{equation}
	Thus, although the number of components depends on $j$, the function $V_j^\vartheta(x)$ is always scalar-valued. We now state the kernel estimates that will be used to prove Proposition \ref{prop:reduction}.
	
	For $j,k\geq1$, $1\leq m\leq d_j$ and $1\leq\ell\leq d_k$, define
	\begin{equation}\label{Ajk}
		A_{jk}^{m,\ell}(u,v,u',v')=\int_{\R^n}\mathcal K_{j,m}(x-u,-v)\overline{\mathcal K_{k,\ell}(x-u',-v')}\,dx.
	\end{equation}

	\begin{lemma}\label{lemma:kernelforL2}
		Let $j,k\geq1$, $1\leq m\leq d_j$ and $1\leq\ell\leq d_k$. Then the following estimates hold uniformly in  $m$ and $\ell$.
		\begin{enumerate}
			\item Suppose $n\geq2$. Then
			\begin{equation}\label{diagonalhighdim}
				|A_{jj}^{m,\ell}(u,v,u',v')|\lesssim 2^{-3nj}(1+|u-u'|)^{-\frac{n-1}{2}}.
			\end{equation}
			and if $k\geq j+4$,
			\begin{equation}\label{offdiagnd}
				|A_{jk}^{m,\ell}(u,v,u',v')|\lesssim 2^{-nj}2^{-2nk}2^{-\frac{j(n-1)}{2}}.
			\end{equation}
			
			\item Suppose $n=1$. If $|v|,|v'|\geq h$ and $0<h\leq2^j$, then
			\begin{equation}\label{diagonal1d}
				|A_{jj}^{m,\ell}(u,v,u',v')|\lesssim 2^{-3j}\left(1+\frac{|u-u'|h^2}{2^{2j}}\right)^{-\frac12}.
			\end{equation}
			If $k\geq j+4$, $|v|\geq h$ and $0<h\leq2^j$, then
			\begin{equation}\label{offdiag1d}
				|A_{jk}^{m,\ell}(u,v,u',v')|\lesssim 2^{-j}2^{-2k}\left(1+\frac{h^2}{2^j}\right)^{-\frac12}.
			\end{equation}
		\end{enumerate}
		In particular, for every $n\geq1$, if $|v|,|v'|\geq h$ and $0<h\leq2^j$, then
		\begin{equation}\label{diagonalkernel}
			|A_{jj}^{m,\ell}(u,v,u',v')|\lesssim 2^{-3nj}
			\left(1+\frac{|u-u'|h^2}{2^{2j}}\right)^{-\frac12}.
		\end{equation}
		Similarly, if $k\geq j+4$, $|v|\geq h$ and $0<h\leq2^j$, then
		\begin{equation}\label{offdiagonalkernel}
			|A_{jk}^{m,\ell}(u,v,u',v')|\lesssim 2^{-nj}2^{-2nk}
			\left(1+\frac{h^2}{2^j}\right)^{-\frac12}.
		\end{equation}
	\end{lemma}
	
	We will prove the Lemma \ref{lemma:kernelforL2} in Section \ref{sec:kernel}. We first show how it implies the desired $L^2$ bound in Proposition \ref{prop:reduction}.
	
	\begin{proof}[Proof of Proposition \ref{prop:reduction}]
		We fix $\vartheta$ and omit it from the notation. In particular, we write $\widetilde F_{j,m}=\widetilde F_{j,m}^{2,\vartheta}$ and $V_j=V_j^\vartheta$. Expanding the $L^2$ norm, we have
		\[\left\|\sum_{j\geq1}V_j\right\|_2^2=\sum_{j,k\geq1}\langle V_j,V_k\rangle.\]
		
		\medskip
		\noindent\textit{\textbf{The $TT^*$ identity.}}
		By \eqref{mainreductioncomponents}, we can write the inner product as
		\begin{align*}
			\langle V_j,V_k\rangle&= \int_{\R^n}V_j(x)\overline{V_k(x)}\,dx\\
			&= \gamma_j\gamma_k\sum_{m=1}^{d_j}\sum_{\ell=1}^{d_k}\iiiint_{\R^{4n}}\widetilde F_{j,m}(u,v)\overline{\widetilde F_{k,\ell}(u',v')}\\
			&\hspace{2 cm}\times\left(\int_{\R^n}\mathcal K_{j,m}((\sqrt2x,0)-(u,v))\overline{\mathcal K_{k,\ell}((\sqrt2x,0)-(u',v'))}\,dx\right)\,du\,dv\,du'\,dv'.
		\end{align*}
		By the change of variables $x\mapsto \frac{x}{\sqrt2}$ in the inner integral, we obtain
		\begin{equation}\label{ttstar}
			\langle V_j,V_k\rangle=2^{-\frac n2}\gamma_j\gamma_k\sum_{m=1}^{d_j}\sum_{\ell=1}^{d_k}\int_{\R^{4n}}A_{jk}^{m,\ell}(u,v,u',v')\widetilde F_{j,m}(u,v)\overline{\widetilde F_{k,\ell}(u',v')}\,du\,dv\,du'\,dv'.
		\end{equation}
		Since the estimates for $|A_{jk}^{m,\ell}|$ in Lemma \ref{lemma:kernelforL2} are uniform in $m,\ell$, and $d_j\leq n$, we can apply the Cauchy–Schwarz inequality for sums over $m,\ell$ in \eqref{ttstar} to obtain
		\begin{equation}\label{componentsnorm}
			\sum_{m=1}^{d_j}|\widetilde F_{j,m}^{2,\vartheta}(u,v)|\leq \sqrt{d_j}\left(\sum_{m=1}^{d_j}|\widetilde F_{j,m}^{2,\vartheta}(u,v)|^2\right)^{\frac{1}{2}}\leq \sqrt n\,|\widetilde F_j^{2,\vartheta}(u,v)|.
		\end{equation}
		Thus, after summation over $m,\ell$ in \eqref{ttstar}, the resulting estimate has the form of a scalar $TT^*$ estimate, with an implicit constant depending only on $n$.
		
		We estimate the near-diagonal terms with $|j-k|<4$ and off-diagonal terms with $|j-k|\geq4$ separately.
		\medskip
		\noindent\textit{\textbf{Near-diagonal terms.}}
		Applying Cauchy--Schwarz in $x$, we have
		\begin{equation}\label{nearindexsum}
			\sum_{\substack{j,k\geq1\\|j-k|<4}}|\langle V_j,V_k\rangle|\leq \sum_{\substack{j,k\geq1\\|j-k|<4}}
			\|V_j\|_2\|V_k\|_2\lesssim \sum_{j\geq1}\|V_j\|_2^2.
		\end{equation}
		
		For a fixed $j$, if the integral in \eqref{Ajk} is nonzero for $k=j$, then the two kernels have supports of size $O(2^j)$. Therefore $|u-u'|\leq|x-u|+|x-u'|\lesssim2^j$. Also, by \eqref{separation}, both $|v|$ and $|v'|$ are at least $h_j=2^{j-\frac{s}{4}}$ on the supports of the modified functions $\widetilde F_j^{2,\vartheta}$. Since $\frac{h_j^2}{2^{2j}}=2^{-\frac{s}{2}}$, the first part of Lemma \ref{lemma:kernelforL2} gives, for every $m,\ell$,
		\begin{equation}\label{Ajjbound}
			|A_{jj}^{m,\ell}(u,v,u',v')|
			\lesssim2^{-3nj}
			\left(1+2^{-\frac{s}{2}}|u-u'|\right)^{-\frac12}
			\mathbf1_{\{|u-u'|\lesssim2^j\}}.
		\end{equation}
		Combining \eqref{ttstar}, \eqref{componentsnorm} and \eqref{Ajjbound}, we obtain
		\begin{equation}\label{diagonalprelim}
			\|V_j\|_2^2\lesssim \gamma_j^2 2^{-3nj}\iint |\widetilde F_j(u',v')|\left(\int_{|u-u'|\lesssim2^j}\int_{\R^n}\frac{|\widetilde F_j(u,v)|}{(1+2^{-\frac{s}{2}}|u-u'|)^{\frac{1}{2}}}\,dv\,du\right)\,du'\,dv'.
		\end{equation}
		
		We claim that, for every $u_0\in\R^n$,
		\begin{equation}\label{diagonalclaim}
			\int_{|u-u_0|\lesssim2^j}\int_{\R^n}\frac{|\widetilde F_j(u,v)|}{(1+2^{-\frac{s}{2}}|u-u_0|)^{\frac{1}{2}}}\,dv\,du\lesssim \lambda2^{2nj}\left[2^{-ns}+\min\left\{1,2^{\frac{s}{4}-\frac j2}\right\}\right].
		\end{equation}
		If we assume \eqref{diagonalclaim} for the moment and substitute it into \eqref{diagonalprelim} along with $\|\widetilde F_j\|_1\lesssim\|F_j\|_1$, we get
		\[\|V_j\|_2^2\lesssim \lambda2^{-nj}\|F_j\|_1\gamma_j^2
		\left[2^{-ns}+\min\left\{1,2^{\frac{s}{4}-\frac j2}\right\}\right].\]
		Note that the term $\gamma_j^22^{-ns}$ is at most $2^{-ns}$. For the other term, we consider three ranges. If $j<\frac{s}{2}$, then $\gamma_j=2^{j-s}$ and
		\[\gamma_j^2\min\left\{1,2^{\frac{s}{4}-\frac j2}\right\}\leq2^{2(j-s)}\leq2^{-s}.\]
		If $\frac{s}{2}\leq j<s$, then
		\[\gamma_j^2 2^{\frac{s}{4}-\frac j2}=2^{\frac{3j}{2}-\frac{7s}{4}}\leq2^{-\frac{s}{4}}.\]
		Finally, if $j\geq s$, then $\gamma_j=1$ and
		\[2^{\frac{s}{4}-\frac j2}\leq2^{-\frac{s}{4}}.\]
		Hence
		\begin{equation}\label{diagsum}
			\sum_{j\geq1}\|V_j\|_2^2\lesssim 2^{-\frac{s}{4}}\lambda\sum_{j\geq1}2^{-nj}\|F_j\|_1\lesssim 2^{-\frac{s}{4}}\lambda^{\frac32},
		\end{equation}
		where the last inequality follows from \eqref{estisumFj}. Thus, it remains to prove the claim \eqref{diagonalclaim}.
		
		\subsection*{Proof of claim \eqref{diagonalclaim}:}
		If $2^j\lesssim2^{\frac{s}{2}}$, then $(1+2^{-\frac{s}{2}}|u-u_0|)^{-\frac{1}{2}}\leq1$. Thus, applying \eqref{localmass} with $r\sim2^j$, we obtain
		\[\int_{|u-u_0|\lesssim2^j}\int\frac{|\widetilde F_j(u,v)|}{(1+2^{-\frac{s}{2}}|u-u_0|)^{\frac{1}{2}}}\,dv\,du\lesssim\lambda2^{2nj}.\]
		This proves the bound corresponding to the term $1$ in the minimum in \eqref{diagonalclaim}.
		
		Now suppose $2^{\frac{s}{2}}\lesssim2^j$. We split the integration region into the ball $|u-u_0|\leq2^{\frac{s}{2}}$ and the annuli
		\[2^{\frac{s}{2}+\ell}<|u-u_0|\leq2^{\frac{s}{2}+\ell+1},\qquad \ell=0,1,\ldots,\quad 2^{\frac{s}{2}+\ell}\lesssim2^j.\]
		Then, $(1+2^{-\frac{s}{2}}|u-u_0|)^{-\frac{1}{2}}=O(2^{-\frac{\ell}{2}})$ on the $\ell$-th annulus. Using \eqref{localmass} on the balls $|u-u_0|\leq2^{\frac{s}{2}+\ell+1}$ containing these annuli,
		we obtain
		\begin{align*}
			\int_{|u-u_0|\lesssim2^j}\int\frac{|\widetilde F_j(u,v)|}{(1+2^{-\frac{s}{2}}|u-u_0|)^{\frac{1}{2}}}\,dv\,du&\lesssim \lambda2^{nj}\left[(2^{\frac{s}{2}}+2^{j-s})^n+\sum_{2^{\frac{s}{2}+\ell}\lesssim2^j}2^{-\frac{\ell}{2}}(2^{\frac{s}{2}+\ell}+2^{j-s})^n\right]\\
			&\lesssim \lambda2^{nj}\left[2^{n(j-s)}+2^{\frac{ns}{2}}\sum_{2^{\frac{s}{2}+\ell}\lesssim2^j}2^{\ell(n-\frac12)}\right]\\
			&\lesssim \lambda2^{nj}\left[2^{n(j-s)}+2^{nj}2^{\frac{s}{4}-\frac j2}\right].
		\end{align*}
		This proves \eqref{diagonalclaim} in the second case.
		
		\medskip
		\noindent\textit{\textbf{Off-diagonal terms.}}
		By symmetry of the inner product, it is enough to consider
		$k\geq j+4$. For fixed $(u',v')$, the integral defining
		$A_{jk}^{m,\ell}(u,v,u',v')$ can be nonzero only if
		$|x-u|\lesssim2^j$ and $|x-u'|\lesssim2^k$ for some
		$x\in\R^n$. Thus
		\[|u-u'|\lesssim2^k.\]
		Also, on the support of $\widetilde F_j$, we have $|v|\lesssim2^j$. From \eqref{offdialemmaesti} for $\widetilde F_j^{2,\vartheta}$ with $u_0=u'$ and $r\sim2^k$, we obtain
		\[\int_{|u-u'|\lesssim2^k}\int_{\R^n}|\widetilde F_j^{2,\vartheta}(u,v)|\,dv\,du\lesssim\lambda2^{nj}2^{nk},\qquad \text{uniformly in }\vartheta.\]
		
		Next, from the second part of Lemma \ref{lemma:kernelforL2} with
		$h=h_j=2^{j-\frac{s}{4}}$, we have
		\[|A_{jk}^{m,\ell}(u,v,u',v')|\lesssim2^{-nj}2^{-2nk}(1+2^{j-\frac{s}{2}})^{-\frac12}.\]
		Combining this bound with \eqref{componentsnorm}, \eqref{offdialemmaesti}, and $\|\widetilde F_k\|_1\lesssim\|F_k\|_1$, we obtain
		\begin{align*}
			|\langle V_j,V_k\rangle|&\lesssim \gamma_j\gamma_k2^{-nj}2^{-2nk}(1+2^{j-\frac{s}{2}})^{-\frac12}\lambda2^{nj}2^{nk}\|F_k\|_1\\
			&\lesssim \lambda\gamma_j(1+2^{j-\frac{s}{2}})^{-\frac12}2^{-nk}\|F_k\|_1,
		\end{align*}
		where we used the fact that $\gamma_k\leq1$.
		
		It remains to sum the factors depending on $j$. Since $\gamma_j=2^{j-s}$ for $j<s$ and $\gamma_j=1$ for $j\geq s$, we obtain
		\begin{align*}
			\sum_{j\geq1}\gamma_j(1+2^{j-\frac{s}{2}})^{-\frac12}&\lesssim \sum_{1\leq j<s}2^{\frac j2-\frac{3s}{4}}+\sum_{j\geq\max\{1,s\}}2^{-\frac j2+\frac s4}\\
			&\lesssim2^{-\frac{s}{4}}.
		\end{align*}
		Hence
		\begin{equation}\label{offdiagestimate}
			\sum_{k\geq1}\sum_{j\leq k-4}|\langle V_j,V_k\rangle|\lesssim2^{-\frac{s}{4}}\lambda\sum_{k\geq1}2^{-nk}\|F_k\|_1\lesssim2^{-\frac{s}{4}}\lambda^{\frac32}.
		\end{equation}
		The terms with $j\geq k+4$ have the same bound by symmetry.
		
		Combining \eqref{nearindexsum}, \eqref{diagsum}, and \eqref{offdiagestimate}, we obtain
		\[\sum_{j,k\geq1}|\langle V_j,V_k\rangle|\lesssim2^{-\frac{s}{4}}\lambda^{\frac32}.\]
		This proves \eqref{propvj} and completes the proof of Proposition \ref{prop:reduction}.
	\end{proof}
	
	\section{Oscillatory kernel estimates}\label{sec:kernel}
	We now prove Lemma \ref{lemma:kernelforL2} by estimating the correlations of the dyadic bilinear Bochner--Riesz kernels. Using the Bessel function expansion, we write these correlations as oscillatory integrals. The main difference from the linear case is that the bilinear kernels are defined on $\R^{2n}$, while their correlations are integrated over $\R^n$. For $n\geq2$, following the ideas of Christ \cite{Christ88} for the linear case, we obtain the required decay using angular stationary phase. The case $n=1$, which does not arise in the linear case, requires a different approach based on van der Corput's lemma.
	
	At the critical index, the order of the Bessel function is $2n-\frac12$. Using the recurrence formulas for Bessel functions (see \cite[Appendix B, Section B.2]{Gra14}), we obtain the finite expansion
	\[K(u,v)=\sum_{\tau=\pm1}e^{2\pi i\tau\sqrt{|u|^2+|v|^2}}\sum_{m=0}^{2n-1}c_{\tau,m}(|u|^2+|v|^2)^{-n-\frac m2},\]
	where $c_{\tau,m}$ are constants depending only on $n,\tau$ and $m$.
	
	Consequently, after multiplying by the dyadic cutoff $\psi(2^{-j}\cdot)$, we can write
	\begin{equation}\label{kernelexpansion}
		K_j(u,v)=2^{-2nj}\sum_{\tau=\pm1}e^{2\pi i\tau\sqrt{|u|^2+|v|^2}}\mathfrak a_{j,\tau}(2^{-j}u,2^{-j}v),
	\end{equation}
	where
	\[\mathfrak a_{j,\tau}(u,v)=\psi(u,v)\sum_{m=0}^{2n-1}c_{\tau,m}2^{-jm}(|u|^2+|v|^2)^{-n-\frac m2}.\]
	Since $\psi$ is smooth and supported in $B(0,2)\setminus B(0,\frac{1}{2})$, the functions $\mathfrak a_{j,\tau}$ are smooth, supported in the same annulus, and satisfy
	\begin{equation}\label{amplitudesymbol}
		\sup_{j\geq1}\|\mathfrak a_{j,\tau}\|_{C^N}\lesssim1
	\end{equation}
	for every nonnegative integer $N$.
	
	Recall that $\mathcal K_{j,m}=K_j$ when $j\geq s$ and $m=1$, whereas for $j<s$, $\mathcal K_{j,m}=\frac1{\sqrt2}(\partial_{u_m}+\partial_{v_m})K_j$. Note that the differentiation of \eqref{kernelexpansion} changes only the amplitude as
	\[(\partial_{u_m}+\partial_{v_m}) e^{2\pi i\tau\sqrt{|u|^2+|v|^2}}=2\pi i\tau\frac{u_m+v_m}{\sqrt{|u|^2+|v|^2}} e^{2\pi i\tau\sqrt{|u|^2+|v|^2}}.\]
	More precisely, define
	\[\mathfrak b_{j,\tau,m}(u,v)=
	\begin{cases}
		\frac{2\pi i\tau(u_m+v_m)}{\sqrt2\sqrt{|u|^2+|v|^2}}\mathfrak a_{j,\tau}(u,v)+\frac{2^{-j}}{\sqrt2} (\partial_{u_m}+\partial_{v_m})\mathfrak a_{j,\tau}(u,v), &j<s,\\
		\mathfrak a_{j,\tau}(u,v),&j\geq s,\quad m=1.
	\end{cases}\]
	Then
	\[\mathcal K_{j,m}(u,v)=2^{-2nj}\sum_{\tau=\pm1} e^{2\pi i\tau\sqrt{|u|^2+|v|^2}} \mathfrak b_{j,\tau,m}(2^{-j}u,2^{-j}v).\]
	Since $\frac{u_m+v_m}{\sqrt{|u|^2+|v|^2}}$ is smooth on the support of $\mathfrak a_{j,\tau}$, $\mathfrak b_{j,\tau,m}$ satisfies the same bounds \eqref{amplitudesymbol}, uniformly in $m$.
	
	For the product of two kernels, write
	\[\mathfrak a_{jk,\tau,\tau'}^{m,\ell}(x,u,v,u',v')=\mathfrak b_{j,\tau,m} (2^{-j}(x-u),-2^{-j}v) \overline{\mathfrak b_{k,\tau',\ell} (2^{-k}(x-u'),-2^{-k}v')}.\]
	Thus, we can rewrite $A_{jk}^{m,\ell}$ in \eqref{Ajk} as
	\begin{equation}\label{Ajkoscillatoryexpansion}
		A_{jk}^{m,\ell}(u,v,u',v')=2^{-2nj}2^{-2nk}\sum_{\tau,\tau'=\pm1}
		\int_{\R^n}
		e^{2\pi i[\tau\sqrt{|x-u|^2+|v|^2}
			-\tau'\sqrt{|x-u'|^2+|v'|^2}]}\mathfrak a_{jk,\tau,\tau'}^{m,\ell}
		(x,u,v,u',v')\,dx.
	\end{equation}
	The phase in the above oscillatory integral is given by
	\begin{equation}\label{phasegeneral}
		\Phi_{\tau,\tau'}(x)=\tau\sqrt{|x-u|^2+|v|^2} -\tau'\sqrt{|x-u'|^2+|v'|^2},
	\end{equation}
	and $\mathfrak a_{jk,\tau,\tau'}^{m,\ell}$ vanishes unless
	\[2^{j-1}\leq\sqrt{|x-u|^2+|v|^2}\leq2\cdot2^j,\qquad 2^{k-1}\leq\sqrt{|x-u'|^2+|v'|^2}\leq2\cdot2^k.\]
	By the chain rule and the bounds for $\mathfrak b_{j,\tau,m}$, we get
	\begin{equation}\label{amplitudederiv}
		\big|\partial_x^\beta\mathfrak a_{jk,\tau,\tau'}^{m,\ell}(x,u,v,u',v')\big|\lesssim 2^{-j|\beta|},\qquad k\geq j.
	\end{equation}

    \medskip
	\begin{proof}[Proof of Lemma \ref{lemma:kernelforL2}]
		We prove the estimates separately for $n\geq2$ and $n=1$.
		
		\medskip
		\subsection{$n\geq2$}
		We first prove the diagonal estimate, followed by the off-diagonal estimate.
		
		\medskip
		\noindent\textit{\textbf{Diagonal estimate $k=j$.}} The support conditions of $\mathfrak a_{jk,\tau,\tau'}^{m,\ell}$ gives $|u-u'|\lesssim2^j$. Moreover, it restricts $x$ to a ball of radius $C2^j$. When $|u-u'|\leq1$, from \eqref{Ajkoscillatoryexpansion}, we get
		\[|A_{jj}^{m,\ell}|\lesssim2^{-3nj}.\]
		Therefore \eqref{diagonalhighdim} holds in this case. We now assume $|u-u'|>1$. By a change of variables into polar coordinates centered at $u'$, we write
		\[x=u'+r\omega,\qquad r>0,\quad\omega\in S^{n-1}.\]
		Then $dx=r^{n-1}\,dr\,d\sigma(\omega)$ and
		\[\Phi_{\tau,\tau'}(u'+r\omega)=\tau\sqrt{r^2+|u-u'|^2-2r\omega\cdot(u-u')+|v|^2}-\tau'\sqrt{r^2+|v'|^2}.\]
		Consequently, \eqref{Ajkoscillatoryexpansion} can be written as
		\begin{align}\label{Ajjoscillatoryexpansionpolar}
			A_{jj}^{m,\ell}(u,v,u',v')&= 2^{-4nj}\sum_{\tau,\tau'=\pm1}\int_0^\infty\int_{S^{n-1}}e^{2\pi i\left[\tau\sqrt{r^2+|u-u'|^2-2r\omega\cdot(u-u')+|v|^2}-\tau'\sqrt{r^2+|v'|^2}\right]}\\
			&\hspace{4cm}\times \mathfrak a_{jj,\tau,\tau'}^{m,\ell}(u'+r\omega,u,v,u',v')\,d\sigma(\omega)\,r^{n-1}\,dr.\nonumber
		\end{align}
		Since the amplitude is supported where $r\lesssim2^j$, the radial integral is taken over $0<r\lesssim2^j$. Let $e=\frac{u-u'}{|u-u'|}$. Differentiating the phase with respect to $\omega$, we obtain
		\[\nabla_{\omega}\Phi_{\tau,\tau'}(u'+r\omega)=-\tau\frac{r|u-u'|}{\sqrt{r^2+|u-u'|^2-2r\omega\cdot(u-u')+|v|^2}}\big(e-(e\cdot\omega)\omega\big).\]
		Thus the critical points are $\omega=\pm e$. Let $U$ be a small neighborhood of these critical points. For $\theta=(\theta_1,\dots,\theta_{n-1})$ and $|\theta|\ll1$, we write $\omega\in U$ as
		\[\omega=\pm e\sqrt{1-|\theta|^2}+\theta_1e_1+\cdots+\theta_{n-1}e_{n-1},\]
		where $e_1,\ldots,e_{n-1}$ form an orthonormal basis of $e^\perp$. Since $\omega\cdot(u-u')=\pm|u-u'|\sqrt{1-|\theta|^2}$, differentiating the phase function twice, we get
		\[D_\theta^2\Phi_{\tau,\tau'}=\pm\tau\frac{r|u-u'|}{R_\pm(\theta)}\bigl(I_{n-1}+O(|\theta|^2)\bigr),\]
		where $I_{n-1}$ is the identity matrix of order $(n-1)\times (n-1)$ and
		\[R_\pm(\theta)=\sqrt{r^2+|u-u'|^2\mp2r|u-u'|\sqrt{1-|\theta|^2}+|v|^2}.\]
		When $\theta=0$, $R_\pm(0)\sim2^j$ on the support of $\mathfrak a_{jk,\tau,\tau'}^{m,\ell}$. Therefore
		\[\left|\det D_\theta^2\Phi_{\tau,\tau'}(u'+r\omega)\big|_{\theta=0}\right|\sim \left(\frac{r|u-u'|}{2^j}\right)^{n-1}.\]
		When $\omega\in U$, $\theta$ is close to 0. Consequently, $R_\pm(\theta)\sim 2^j$ and $r|u-u'|\lesssim2^{2j}$ on the support of $\mathfrak a_{jk,\tau,\tau'}^{m,\ell}$. Thus
		\[\left|\det D_\theta^2\Phi_{\tau,\tau'}(u'+r\omega)\right|\gtrsim\left(\frac{r|u-u'|}{2^j}\right)^{n-1}.\]
		
		Observe that from \eqref{amplitudederiv}, we have
		\[\big|\partial_\omega^\beta \mathfrak a_{jj,\tau,\tau'}^{m,\ell}(u'+r\omega,u,v,u',v')\big|\lesssim 1,\qquad r\lesssim2^j.\]
		Also, for every $|\beta|\geq1$, we have
		\[\big|\partial_\omega^\beta \Phi_{\tau,\tau'}(u'+r\omega)\big|\lesssim \frac{r|u-u'|}{2^j}.\]
		Thus, when $\frac{r|u-u'|}{2^j}\geq1$, by stationary-phase method, we obtain
		\begin{align*}
			&\left|\int_{U}e^{2\pi i[\tau\sqrt{r^2+|u-u'|^2-2r\omega\cdot(u-u')+|v|^2}-\tau'\sqrt{r^2+|v'|^2}]}\mathfrak a_{jj,\tau,\tau'}^{m,\ell}(u'+r\omega,u,v,u',v')\,d\sigma(\omega)\right|\\
			&\lesssim \left(1+\frac{r|u-u'|}{2^j}\right)^{-\frac{n-1}{2}}.
		\end{align*}
		
		Away from  the critical points on $U^{c}$, we have
		\[\big|\nabla_{\omega}\Phi_{\tau,\tau'}(u'+r\omega)\big|\gtrsim\frac{r|u-u'|}{2^j}.\]
		Therefore, when $\frac{r|u-u'|}{2^j}\geq1$, by an integration by parts argument, we obtain
		\begin{align*}
			&\left|\int_{U^{c}}e^{2\pi i[\tau\sqrt{r^2+|u-u'|^2-2r\omega\cdot(u-u')+|v|^2}-\tau'\sqrt{r^2+|v'|^2}]}\mathfrak a_{jj,\tau,\tau'}^{m,\ell}(u'+r\omega,u,v,u',v')\,d\sigma(\omega)\right|\\
			&\lesssim \left(1+\frac{r|u-u'|}{2^j}\right)^{-\frac{n-1}{2}}.
		\end{align*}
		When $\frac{r|u-u'|}{2^j}<1$, the above estimate follows directly by integrating \eqref{Ajkoscillatoryexpansion}.
		Combining the estimates for $\frac{r|u-u'|}{2^j}<1$ and $\frac{r|u-u'|}{2^j}\geq1$ including both near and away from critical points, we obtain
		\begin{align}\label{angularbound}
			&\left|\int_{S^{n-1}} e^{2\pi i[\tau\sqrt{r^2+|u-u'|^2-2r\omega\cdot(u-u')+|v|^2} -\tau'\sqrt{r^2+|v'|^2}]} \mathfrak a_{jj,\tau,\tau'}^{m,\ell}(u'+r\omega,u,v,u',v')\,d\sigma(\omega)\right|\\
			&\lesssim \left(1+\frac{r|u-u'|}{2^j}\right)^{-\frac{n-1}{2}}.\nonumber
		\end{align}
		
		Substituting \eqref{angularbound} into \eqref{Ajjoscillatoryexpansionpolar} and then splitting the integral in $r$ at $r=\frac{2^j}{|u-u'|}$, we obtain
		\begin{align*}
			|A_{jj}^{m,\ell}(u,v,u',v')|&\lesssim 2^{-4nj}\int_0^{C2^j}r^{n-1}\left(1+\frac{r|u-u'|}{2^j}\right)^{-\frac{n-1}{2}}\,dr\\
			&\lesssim2^{-4nj}\left[\left(\frac{2^j}{|u-u'|}\right)^n+\left(\frac{2^j}{|u-u'|}\right)^{\frac{n-1}{2}}\int_0^{C2^j}r^{\frac{n-1}{2}}\,dr\right]\\
			&\lesssim2^{-3nj}|u-u'|^{-\frac{n-1}{2}},
		\end{align*}
		where we used the fact that $|u-u'|>1$. This completes the proof of \eqref{diagonalhighdim}. Since $n\geq2$ and $0<h\leq2^j$, we have
		\[(1+|u-u'|)^{-\frac{n-1}{2}}\leq\left(1+\frac{|u-u'|h^2}{2^{2j}}\right)^{-\frac12},\]
		which also proves \eqref{diagonalkernel} in this range.
		
		\medskip
		\noindent\textit{\textbf{Off-diagonal estimate $k\geq j+4$.}}
		As in the diagonal case, we use polar coordinates centered at $u'$ and rewrite \eqref{Ajkoscillatoryexpansion} as
		\begin{align}\label{Ajkoscillatoryexpansionpolar}
			&A_{jk}^{m,\ell}(u,v,u',v')\\
			&=2^{-2nj}2^{-2nk}\sum_{\tau,\tau'=\pm1}\int_0^\infty\int_{S^{n-1}} e^{2\pi i\Phi_{\tau,\tau'}(u'+r\omega)}\mathfrak a_{jk,\tau,\tau'}^{m,\ell}
			(u'+r\omega,u,v,u',v')\,d\sigma(\omega)\,r^{n-1}\,dr.\nonumber
		\end{align}
		From the support of $\mathfrak a_{jk,\tau,\tau'}^{m,\ell}$, we have
		\[|r-|u-u'||\leq2^{j+1},\qquad \sqrt{r^2+|v'|^2}\geq2^{k-1}.\]
		We consider the integrals in two cases.
		
		\smallskip
		\noindent\textit{Case 1: $|u-u'|\leq2^{j+2}$.}
		In this case, we have $r\leq|u-u'|+2^{j+1}\leq6\cdot2^j$.
		Hence, we can obtain the following estimate similar to \eqref{angularbound}.
		\[\left|\int_{S^{n-1}}e^{2\pi i\Phi_{\tau,\tau'}(u'+r\omega)}\mathfrak a_{jk,\tau,\tau'}^{m,\ell}(u'+r\omega,u,v,u',v')\,d\sigma(\omega)\right| \lesssim\left(1+\frac{r|u-u'|}{2^j}\right)^{-\frac{n-1}{2}}.\]
		With this estimate, it remains to compute the integral in $r$. We first consider
		\[r\leq 4\max\{|u-u'|,2^{\frac{j}{2}}\}.\]
		If $|u-u'|\leq2^{\frac{j}{2}}$, then we have from the trivial bound,
		\[\int_0^{C2^{\frac{j}{2}}}r^{n-1}\,dr\lesssim2^{\frac{jn}{2}}\]
		If $2^{\frac{j}{2}}<|u-u'|\leq2^{j+2}$, splitting the integral at $r=\frac{2^j}{|u-u'|}$, we obtain
		\begin{align*}
			\int_0^{C|u-u'|}r^{n-1}\left(1+\frac{r|u-u'|}{2^j}\right)^{-\frac{n-1}{2}}\,dr&\lesssim \left(\frac{2^j}{|u-u'|}\right)^n+\left(\frac{2^j}{|u-u'|}\right)^{\frac{n-1}{2}}\int_0^{C|u-u'|}r^{\frac{n-1}{2}}\,dr\\
			&\lesssim \frac{2^{nj}}{|u-u'|^n}+2^{\frac{j(n-1)}2}|u-u'|\\
			&\lesssim2^{\frac{j(n+1)}2}.
		\end{align*}
		
		Next, we consider
		\[r\geq 4\max\{|u-u'|,2^{\frac{j}{2}}\}.\]
		From the support of $\mathfrak a_{jk,\tau,\tau'}^{m,\ell}$, we obtain
		\begin{align*}
			\tau\partial_r\Phi_{\tau,\tau'}(u'+r\omega)&= \frac{r-\omega\cdot(u-u')}{\sqrt{|r\omega-(u-u')|^2+|v|^2}}-\tau\tau'\frac{r}{\sqrt{r^2+|v'|^2}}\\
			&\geq \frac{r-|u-u'|}{2^{j+1}}-\frac{r}{2^{k-1}}\\
			&\geq \frac{r}{4\cdot2^j},
		\end{align*}
		where we used the fact that $k\geq j+4$. Moreover, we have the upper bounds
		\[|\partial_r\Phi_{\tau,\tau'}(u'+r\omega)|\lesssim\frac{r}{2^j},\]
		and for every integer $M\geq2$,
		\[|\partial_r^M\Phi_{\tau,\tau'}(u'+r\omega)|\lesssim 2^{j(1-M)}.\]
		Since $r\lesssim2^j$, differentiating $\frac1{\partial_r\Phi_{\tau,\tau'}(u'+r\omega)}$, we obtain
		\[\left|\partial_r^M\left(\frac1{\partial_r\Phi_{\tau,\tau'}(u'+r\omega)}\right)\right|\lesssim \frac{2^j}{r^{M+1}},\qquad M\geq0.\]
		Also, by \eqref{amplitudederiv}, we have
		\[\left|\partial_r^M\mathfrak a_{jk,\tau,\tau'}^{m,\ell}(u'+r\omega,u,v,u',v')\right|\lesssim 2^{-jM}\lesssim r^{-M}.\]
		
		Since $e^{2\pi i\Phi_{\tau,\tau'}}=\frac1{2\pi i\partial_r\Phi_{\tau,\tau'}}\partial_r e^{2\pi i\Phi_{\tau,\tau'}}$, integrating by parts in $r$ along with the above derivative estimates, we obtain
		\begin{align*}
			&\left|\int_{r\gtrsim\max\{|u-u'|,2^{\frac{j}{2}}\}}e^{2\pi i\Phi_{\tau,\tau'}(u'+r\omega)}\mathfrak a_{jk,\tau,\tau'}^{m,\ell}(u'+r\omega,u,v,u',v')r^{n-1}\,dr\right|\\
			&\lesssim 2^{jN}\int_{r\gtrsim\max\{|u-u'|,2^{\frac{j}{2}}\}}r^{n-1-2N}\,dr.
		\end{align*}
		We choose $N>\frac{n}{2}$ and obtain
		\begin{align*}
			2^{jN}\int_{r\gtrsim\max\{|u-u'|,2^{\frac{j}{2}}\}}r^{n-1-2N}\,dr&\lesssim 2^{jN}\bigl(\max\{|u-u'|,2^{\frac{j}{2}}\}\bigr)^{n-2N}\\
			&\lesssim 2^{\frac{jn}{2}}.
		\end{align*}
		Since the above bound is uniform in $\omega$, integrating over $S^{n-1}$ we obtain the same bound.
		
		\smallskip
		\noindent\textit{Case 2: $|u-u'|>2^{j+2}$.}
		Let $e=\frac{u-u'}{|u-u'|}$. From the support of $\mathfrak a_{jk,\tau,\tau'}^{m,\ell}$, we have
		\[|r\omega-(u-u')|^2=(r-|u-u'|)^2+2r|u-u'|(1-\omega\cdot e)\lesssim2^{2j}.\]
		Thus, $r\sim|u-u'|$ and
		\[1-\omega\cdot e\lesssim\frac{2^{2j}}{r|u-u'|}.\]
		Since $|\omega|=|e|=1$, this implies
		\[|\omega-(\omega\cdot e)e|\lesssim\frac{2^j}{\sqrt{r|u-u'|}}.\]
		Hence $\omega$ lies in a small neighborhood of $e$. As in the diagonal case, we choose an orthonormal basis $e_1,\ldots,e_{n-1}$ of $e^\perp$. Rescaling the angular coordinates according to the size of this neighborhood, we write
		\[\omega=e\sqrt{1-\frac{2^{2j}}{r|u-u'|}|\theta|^2}+\frac{2^j}{\sqrt{r|u-u'|}}(\theta_1e_1+\cdots+\theta_{n-1}e_{n-1}),\qquad |\theta|\lesssim1.\]
		By change of variables in above coordinates, we get
		\[d\sigma(\omega)=\left(\frac{2^{2j}}{r|u-u'|}\right)^{\frac{n-1}{2}}\left(1-\frac{2^{2j}}{r|u-u'|}|\theta|^2\right)^{-\frac12}d\theta.\]
		Moreover, in these coordinates, the phase function becomes
		\[\Phi_{\tau,\tau'}(u'+r\omega)=\tau\Bigg[(r-|u-u'|)^2+|v|^2+2r|u-u'|\left(1-\sqrt{1-\frac{2^{2j}}{r|u-u'|}|\theta|^2}\right)\Bigg]^{\frac12}-\tau'\sqrt{r^2+|v'|^2}.\]
		Differentiating with respect to $\theta$, we obtain
		\[\nabla_\theta\bigl(2^{-j}\Phi_{\tau,\tau'}(u'+r\omega)\bigr)=\frac{\tau\theta}{2^{-j}\sqrt{|r\omega-(u-u')|^2+|v|^2}\sqrt{1-\frac{2^{2j}}{r|u-u'|}|\theta|^2}}.\]
		Since the first square root is comparable to $2^j$ on the support of the amplitude, the gradient vanishes only when $\theta=0$. 
		
		Let $U$ be a fixed small neighborhood of $\theta=0$. Differentiating once more on the set $U$, we get
		\[D_\theta^2\bigl(2^{-j}\Phi_{\tau,\tau'}(u'+r\omega)\bigr)=\frac{\tau}{2^{-j}\sqrt{|r\omega-(u-u')|^2+|v|^2}}\bigl(I_{n-1}+O(|\theta|^2)\bigr).\]
		Thus, the Hessian is uniformly nondegenerate in $U$. On $U^{c}$, we have
		\[\left|\nabla_\theta\bigl(2^{-j}\Phi_{\tau,\tau'}(u'+r\omega)\bigr)\right|\gtrsim|\theta|\geq\delta.\]
		
		Since $r\sim|u-u'|$ and $\frac{2^j}{|u-u'|}\leq\frac{1}{4}$, for every $|\beta|\geq1$, we have
		\[\left|\partial_\theta^\beta\left(\frac{r\omega-(u-u')}{2^j}\right)\right|\lesssim \sqrt{\frac r{|u-u'|}}\lesssim1.\]
		Similarly, using $r\lesssim2^k$, we have
		\[\left|\partial_\theta^\beta(2^{-k}r\omega)\right|\lesssim1.\]
		Consequently, by the chain rule, for every multi-index $\beta$, we have
		\[\left|\partial_\theta^\beta\mathfrak a_{jk,\tau,\tau'}^{m,\ell}(u'+r\omega,u,v,u',v')\right|+\left|\partial_\theta^\beta\left(1-\frac{2^{2j}}{r|u-u'|}|\theta|^2\right)^{-\frac12}\right|\lesssim1.\]
		
		Therefore, applying stationary phase on $U$ and integration by parts on $U^{c}$, as in the diagonal case, we obtain
		\[\left|\int_{S^{n-1}}e^{2\pi i\Phi_{\tau,\tau'}(u'+r\omega)}\mathfrak a_{jk,\tau,\tau'}^{m,\ell}(u'+r\omega,u,v,u',v')\,d\sigma(\omega)\right|\lesssim \left(\frac{2^{2j}}{r|u-u'|}\right)^{\frac{n-1}{2}}2^{-\frac{j(n-1)}2}.\]
		Since $|r-|u-u'||\leq2^{j+1}$ and
		$r\sim|u-u'|$, integrating in $r$, we get
		\begin{align*}
			\int_{|r-|u-u'||\leq2^{j+1}}r^{n-1}\left(\frac{2^{2j}}{r|u-u'|}\right)^{\frac{n-1}{2}}2^{-\frac{j(n-1)}2}\,dr&= 2^{\frac{j(n-1)}2}\int_{|r-|u-u'||\leq2^{j+1}}\left(\frac r{|u-u'|}\right)^{\frac{n-1}{2}}\,dr\\
			&\lesssim 2^{\frac{j(n+1)}2}.
		\end{align*}
		Combining both cases, each oscillatory integral in \eqref{Ajkoscillatoryexpansion} is bounded by $C2^{\frac{j(n+1)}2}$. Consequently, we obtain
		\begin{align*}
			|A_{jk}^{m,\ell}(u,v,u',v')|
			&\lesssim2^{-2nj}2^{-2nk}2^{\frac{j(n+1)}2}\\
			&=2^{-nj}2^{-2nk}2^{-\frac{j(n-1)}2}.
		\end{align*}
		This completes the proof of \eqref{offdiagnd}. Now, for $n\geq2$ and $0<h\leq2^j$, we have
		\[2^{-\frac{j(n-1)}2}\lesssim\left(1+\frac{h^2}{2^j}\right)^{-\frac12},\]
		which proves \eqref{offdiagonalkernel}.
		
		\subsection{$n=1$}
		In this case, \eqref{Ajkoscillatoryexpansion} is
		\begin{equation}\label{Ajkoscillatoryexpansion1}
			A_{jk}^{m,\ell}(u,v,u',v')=2^{-2j}2^{-2k}\sum_{\tau,\tau'=\pm1}\int_{\R}e^{2\pi i\Phi_{\tau,\tau'}(x)}\mathfrak a_{jk,\tau,\tau'}^{m,\ell}(x,u,v,u',v')\,dx,
		\end{equation}
		where
		\[\Phi_{\tau,\tau'}(x)=\tau\sqrt{(x-u)^2+v^2}-\tau'\sqrt{(x-u')^2+(v')^2}.\]
		Differentiating the phase function, we obtain
		\begin{align*}
			\Phi_{\tau,\tau'}'(x)&=\tau\frac{x-u}{\sqrt{(x-u)^2+v^2}}-\tau'\frac{x-u'}{\sqrt{(x-u')^2+(v')^2}},\\
			\Phi_{\tau,\tau'}''(x)&=\tau\frac{v^2}{((x-u)^2+v^2)^{\frac{3}{2}}}-\tau'\frac{(v')^2}{((x-u')^2+(v')^2)^{\frac{3}{2}}}.
		\end{align*}
		When $\tau=\tau'$, we write
		\[\Phi(x)=\sqrt{(x-u)^2+v^2}-\sqrt{(x-u')^2+(v')^2},\]
		so that $\Phi_{\tau,\tau}=\tau\Phi$.
		
		For fixed $u,v\in\R$ and $j,k\geq1$, $\mathfrak a_{jk,\tau,\tau'}^{m,\ell}$ is supported in $x$ on the region
		\[I_{jk}=\{x:2^{j-1}\leq \sqrt{(x-u)^2+v^2}\leq2^{j+1},\, 2^{k-1}\leq \sqrt{(x-u')^2+(v')^2}\leq2^{k+1}\}.\]
		Each shell condition gives at most two intervals. Thus, $I_{jk}$ is a union of at most four intervals, with total length $O(2^j)$. Note that the support of $\mathfrak a_{jk,\tau,\tau'}^{m,\ell}$ itself need not have finitely many connected components.
		
		On $I_{jk}$, the two square roots are comparable to $2^j$ and $2^k$, respectively. Moreover, by \eqref{amplitudederiv}, we have
		\[\|\mathfrak a_{jk,\tau,\tau'}^{m,\ell}\|_\infty+2^j\|\partial_x\mathfrak a_{jk,\tau,\tau'}^{m,\ell}\|_\infty\lesssim1,\qquad k\geq j.\]
		We apply the first and second derivative tests of van der Corput's lemma.
		
		\medskip
		\noindent\textit{\textbf{Diagonal estimate $k=j$.}}
		Suppose $|v|,|v'|\geq h$. The support conditions of $\mathfrak a_{jk,\tau,\tau'}^{m,\ell}$ give
		\[\sqrt{(x-u)^2+v^2}\sim2^j,\qquad\sqrt{(x-u')^2+(v')^2}\sim2^j,\]
		and $|u-u'|\lesssim2^j$. When $u=u'$, the trivial estimate $|A_{jj}^{m,\ell}|\lesssim2^{-3j}$ proves \eqref{diagonal1d}. Hence we assume $u\ne u'$.
		
		First, suppose $\tau=-\tau'$. Then
		\begin{align*}
			|\Phi_{\tau,\tau'}''(x)|&=\frac{v^2}{((x-u)^2+v^2)^{\frac{3}{2}}}+\frac{(v')^2}{((x-u')^2+(v')^2)^{\frac{3}{2}}}\\
			&\gtrsim\frac{h^2}{2^{3j}}.
		\end{align*}
		Since the second derivative has a fixed sign, by van der Corput's lemma, we get
		\begin{align*}
			\left|\int_\R e^{2\pi i\Phi_{\tau,\tau'}(x)}\mathfrak a_{jj,\tau,\tau'}^{m,\ell}(x,u,v,u',v')\,dx\right|&\lesssim \min\left\{2^j,\left(\frac{2^{3j}}{h^2}\right)^{\frac{1}{2}}\right\}\\
			&\lesssim 2^j\left(1+\frac{h^2}{2^j}\right)^{-\frac12}\\
			&\lesssim 2^j\left(1+\frac{|u-u'|h^2}{2^{2j}}\right)^{-\frac12},
		\end{align*}
		where we used $|u-u'|\lesssim2^j$ in the last inequality.
		
		Next, suppose $\tau=\tau'$. Observe that
		\begin{equation}\label{eq:inversionidentity}
			\frac{x-u}{\sqrt{(x-u)^2+v^2}}\frac{1-\frac{(x-u)^2}{(x-u)^2+v^2}}{\frac{v^2}{((x-u)^2+v^2)^{\frac{3}{2}}}}=x-u.
		\end{equation}
		The same identity holds with $(u,v)$ replaced by $(u',v')$. Since $r(1-r^2)$ is Lipschitz on
		$[-1,1]$ and
		\[\frac{v^2}{((x-u)^2+v^2)^{\frac{3}{2}}},\,\frac{(v')^2}{((x-u')^2+(v')^2)^{\frac{3}{2}}}\gtrsim\frac{h^2}{2^{3j}},\]
		from \eqref{eq:inversionidentity}, we get
		\begin{align*}
			|u-u'|&\leq |x-u-(x-u')|\\
			&\lesssim\frac{2^{3j}}{h^2}\left(|\Phi'(x)|+|x-u'|\,|\Phi''(x)|\right)\\
			&\lesssim\frac{2^{3j}}{h^2}\left(|\Phi'(x)|+2^j|\Phi''(x)|\right).
		\end{align*}
		Thus
		\begin{equation}\label{eq:1djet}
			|\Phi'(x)|+2^j|\Phi''(x)|\gtrsim\frac{|u-u'|h^2}{2^{3j}}.
		\end{equation}
		
		Next, the zeros of $\Phi''$ satisfy the equation
		\[|v|^{\frac43}((x-u')^2+(v')^2)=|v'|^{\frac43}((x-u)^2+v^2).\]
		This is a quadratic equation in $x$. Since $u\ne u'$, it cannot hold identically. Thus, $\Phi''$ has at most two zeros. Dividing the intervals of $I_{jj}$ at these zeros, we obtain a bounded number of intervals on which $\Phi'$ is monotone. Moreover, the support conditions and \eqref{eq:1djet} hold on each of these intervals.
		
		If $|u-u'|h^2\leq2^{2j}$, the trivial bound $|A_{jj}^{m,\ell}|\lesssim2^{-3j}$ is sufficient. Otherwise, we divide each of these intervals into the regions where
		\[|\Phi'(x)|\geq c\frac{|u-u'|h^2}{2^{3j}}\qquad\text{and}\qquad |\Phi'(x)|< c\frac{|u-u'|h^2}{2^{3j}},\]
		where $c>0$ is sufficiently small.
		
		On the first set, by the first derivative test we get
		\[\left|\int e^{2\pi i\tau\Phi(x)}\mathfrak a_{jj,\tau,\tau}^{m,\ell}(x,u,v,u',v')\,dx\right|\lesssim\frac{2^{3j}}{|u-u'|h^2}.\]
		On the second set, from \eqref{eq:1djet}, we have
		\[|\Phi''(x)|\gtrsim\frac{|u-u'|h^2}{2^{4j}}.\]
		Thus, by the second derivative test, we get
		\[\left|\int e^{2\pi i\tau\Phi(x)}\mathfrak a_{jj,\tau,\tau}^{m,\ell}(x,u,v,u',v')\,dx\right|\lesssim \left(\frac{2^{4j}}{|u-u'|h^2}\right)^{\frac{1}{2}}.\]
		Since $|u-u'|h^2>2^{2j}$, the second bound is larger than the first. Combining these estimates, we obtain
		\[\left|\int_\R e^{2\pi i\tau\Phi(x)}\mathfrak a_{jj,\tau,\tau}^{m,\ell}(x,u,v,u',v')\,dx\right|\lesssim 2^j\left(1+\frac{|u-u'|h^2}{2^{2j}}\right)^{-\frac12}.\]
		Substituting the estimates for both sign choices into \eqref{Ajkoscillatoryexpansion1} proves
		\eqref{diagonal1d}.
		
		\medskip
		\noindent\textit{\textbf{Off-diagonal estimate $k\geq j+4$.}}
		Suppose $|v|\geq h$. We follow the same argument as in the diagonal case.
		
		First, suppose $\tau=-\tau'$. Then
		\begin{align*}
			|\Phi_{\tau,\tau'}''(x)|&=\frac{v^2}{((x-u)^2+v^2)^{\frac{3}{2}}}+\frac{(v')^2}{((x-u')^2+(v')^2)^{\frac{3}{2}}}\\
			&\gtrsim\frac{h^2}{2^{3j}}.
		\end{align*}
		Thus, by the second derivative test and the trivial bound, we get
		\[\left|\int_\R e^{2\pi i\Phi_{\tau,\tau'}(x)}\mathfrak a_{jk,\tau,\tau'}^{m,\ell}(x,u,v,u',v')\,dx\right|\lesssim 2^j\left(1+\frac{h^2}{2^j}\right)^{-\frac12}.\]
		
		Now suppose $\tau=\tau'$. Observe that
		\[\frac{(v')^2}{((x-u')^2+(v')^2)^{\frac{3}{2}}}=\frac1{\sqrt{(x-u')^2+(v')^2}}\left(1-\frac{(x-u')^2}{(x-u')^2+(v')^2}\right).\]
		Using $|r^2-(r')^2|\leq2|r-r'|$ for
		$|r|,|r'|\leq1$ and the support conditions of $\mathfrak a_{jk,\tau,\tau'}^{m,\ell}$,
		we obtain
		\begin{align*}
			\frac{(v')^2}{((x-u')^2+(v')^2)^{\frac{3}{2}}}&\leq\frac2{2^k}\left(1-\frac{(x-u)^2}{(x-u)^2+v^2}+2|\Phi'(x)|\right)\\
			&\leq 4\cdot2^{j-k}\frac{v^2}{((x-u)^2+v^2)^{\frac{3}{2}}}+\frac4{2^k}|\Phi'(x)|\\
			&\leq \frac14\frac{v^2}{((x-u)^2+v^2)^{\frac{3}{2}}}+\frac1{4\cdot2^j}|\Phi'(x)|,
		\end{align*}
		where we used $k\geq j+4$ in the last inequality. Consequently, we get
		\[\Phi''(x)\geq \frac34\frac{v^2}{((x-u)^2+v^2)^{\frac{3}{2}}}-\frac1{4\cdot2^j}|\Phi'(x)|.\]
		Since $\frac{v^2}{((x-u)^2+v^2)^{\frac{3}{2}}}\gtrsim\frac{h^2}{2^{3j}}$, we get
		\begin{equation}\label{eq:1doffjet}
			|\Phi'(x)|+2^j|\Phi''(x)|\gtrsim\frac{h^2}{2^{2j}}.
		\end{equation}
		
		As in the diagonal case, the equation $\Phi''(x)=0$ is quadratic in $x$ when $v'\ne0$. It cannot hold identically, since this would imply $u=u'$ and $|v|=|v'|$, which is not possible for $k\geq j+4$. If $v'=0$, then $\Phi''(x)>0$. As in the diagonal case, we divide the intervals of $I_{jk}$ at the zeros of $\Phi''$. This gives a bounded number of intervals on which $\Phi'$ is monotone. Moreover, the support conditions and \eqref{eq:1doffjet} hold on each of these intervals.
		
		If $h^2\leq2^j$, the trivial bound $|A_{jk}^{m,\ell}|\lesssim2^{-j}2^{-2k}$ is enough. Otherwise, on each monotonic interval, we consider the sets where
		\[|\Phi'(x)|\geq c\frac{h^2}{2^{2j}}\qquad\text{and}\qquad |\Phi'(x)|<c\frac{h^2}{2^{2j}}.\]
		On the first set, by the first derivative test we get
		\[\left|\int e^{2\pi i\tau\Phi(x)}\mathfrak a_{jk,\tau,\tau}^{m,\ell}(x,u,v,u',v')\,dx\right|\lesssim \frac{2^{2j}}{h^2}.\]
		On the second set, \eqref{eq:1doffjet} implies
		\[|\Phi''(x)|\gtrsim\frac{h^2}{2^{3j}},\]
		and by the second derivative test we get
		\[\left|\int e^{2\pi i\tau\Phi(x)}\mathfrak a_{jk,\tau,\tau}^{m,\ell}(x,u,v,u',v')\,dx\right|\lesssim\left(\frac{2^{3j}}{h^2}\right)^{\frac{1}{2}}.\]
		Since $h^2>2^j$, the second estimate is larger than the first. Combining both estimates, we obtain
		\[\left|\int_\R e^{2\pi i\tau\Phi(x)}\mathfrak a_{jk,\tau,\tau}^{m,\ell}(x,u,v,u',v')\,dx\right|\lesssim 2^j\left(1+\frac{h^2}{2^j}\right)^{-\frac12}.\]
		Substituting the estimates for both sign choices into \eqref{Ajkoscillatoryexpansion1} proves \eqref{offdiag1d}.
		
		All the estimates are uniform in $m,\ell$ and $\tau,\tau'$. This completes the proof of Lemma \ref{lemma:kernelforL2}.
	\end{proof}
	
	\subsection*{Acknowledgements.} Surjeet Singh Choudhary is supported by the National Center for Theoretical Sciences through National Science and Technology Council of Taiwan under research grant number 115-2124-M-002-009-. Saurabh Shrivastava acknowledges the support from Anusandhan National Research Foundation (ANRF), India under the projects ANRF/ARG/2025/000940/MS. 
	
	\bibliography{biblio}
	
\end{document}